\documentclass{article}
\usepackage{hyperref}
\usepackage[american]{babel}
\usepackage{amsfonts,amsmath,amssymb,epsf,epsfig, stmaryrd}

\usepackage{amsthm}

\usepackage{xypic}
\xyoption{all}
\input{xypic}
\newdir{ >}{{}*!/-8pt/\dir{>}}

\renewcommand{\hat}{\widehat}
\newcommand{\im}{\mathop{{\rm im}}\nolimits}
\newcommand{\Hom}{\mathop{{\rm Hom}}\nolimits}

\newcommand{\coker}{\mathop{{\rm coker}}\nolimits}

\def\N{\mathbb{N}}
\def\C{\mathbb{C}}

     \def\Gl{{\rm Gl}}
   \def\dim{{\rm dim}}
     \def\Hom{{\rm Hom}}

\newtheorem{theo}{Theorem}[section]
\newtheorem{defi}[theo]{Definition}
\newtheorem{rema}[theo]{Remark}
\newtheorem{prop}[theo]{Proposition}
\newtheorem{cor}[theo]{Corollary}

\newtheorem{pro}[theo]{Problem}

\newenvironment{rem}{\begin{rema}\rm}{\end{rema}}

\newenvironment{prf}{\begin{proof}}{\end{proof}}

\begin{document}
\title{A rigid Leibniz algebra with non-trivial $HL^2$}

\author{Bakhrom Omirov\\
National University of Uzbekistan \\
\and Friedrich Wagemann\\
     Universit\'e de Nantes }

\maketitle

\begin{abstract}
In this article, we generalize Richardson's example of a rigid Lie algebra with non-trivial
$H^2$ to the Leibniz setting. Namely, we consider the hemisemidirect product
${\mathfrak h}$ of a semidirect product Lie algebra $M_k\rtimes{\mathfrak g}$
of a simple Lie algebra ${\mathfrak g}$ with some non-trivial irreducible
${\mathfrak g}$-module $M_k$ with a non-trivial irreducible ${\mathfrak g}$-module $I_l$.
Then for ${\mathfrak g}={\mathfrak s}{\mathfrak l}_2(\C)$, we take $M_k$
(resp. $I_l$) to be the standard irreducible ${\mathfrak s}{\mathfrak l}_2(\C)$-module
of dimension $k+1$ (resp. $l+1$). Assume
$\frac{k}{2}>5$ is an odd integer and $l>2$ is odd, then we show that the Leibniz algebra
${\mathfrak h}$ is geometrically rigid and has non-trivial $HL^2$ with adjoint coefficients.
We close the article with an appendix where we record further results on the question whether $H^2({\mathfrak g},{\mathfrak g})=0$ implies $HL^2({\mathfrak g},{\mathfrak g})=0$. 
\end{abstract}

\section*{Introduction}
Let $k$ be a field. In order to study the variety of all $n$-dimensional Lie algebras over
$k$, one fixes a basis $(e_i)_{i=1,\ldots,n}$ in $k^n$ and represents a Lie algebra by its
structure constants $(c_{ij}^k)_{i,j,k\in\{1,\ldots,n\}}$ given by
$$[e_i,e_j]\,=\,\sum_{i=1}^n c_{ij}^ke_k.$$
These structure constants are elements of the vector space $\Hom(\Lambda^2(k^n),k^n)$
which must satisfy the quadratic equations
$$\sum_{p=1}^n(c_{jl}^pc_{ip}^k-c_{ij}^pc_{pl}^k-c_{il}^pc_{jp}^k)\,=\,0$$
for $i,j,l,k\in\{1,\ldots,n\}$ owing to the Jacobi identity.
The group $\Gl(k^n)$ acts on the algebraic variety of structure constants
by base changes. The {\it variety of Lie algebra laws} over $k$ is by definition the quotient
of the algebraic variety defined by the above quadratic equations by the action of $\Gl(k^n)$.
This action is usually badly behaved and the quotient has singularities, is non-reduced etc.
even if $k$ is algebraically closed and of characteristic zero, which we suppose from now on.
Therefore it is not a variety in the usual sense. We nevertheless continue to call it
the variety of Lie algebra laws.

The analoguous picture for Leibniz algebra laws has been explored in \cite{Bal}.
The main change here is that the structure constants are not supposed to be antisymmetric
anymore, thus they lie in $\Hom((k^n)^{\otimes 2},k^n)$.
Balavoine shows that a point in the variety of Leibniz algebra laws is reduced and
geometrically rigid (i.e. the $\Gl(k^n)$-orbit of the corresponding Leibniz algebra
${\mathfrak h}$ is open in the Zariski topology) if and only if its second adjoint
Leibniz cohomology space $HL^2({\mathfrak h},{\mathfrak h})$ is zero (see \cite{Bal}).
The goal of the present paper is to give an example
of a non-reduced point of the variety of Leibniz algebra laws, i.e. we give an example
of a finite-dimensional Leibniz algebra ${\mathfrak h}$ over $k$ which is geometrically rigid,
but which has $HL^2({\mathfrak h},{\mathfrak h})\not=0$. We will use the word {\it rigid} in this article
for geometrically rigid and express algebraic rigidity (in the Lie or Leibniz context) by stating that the second cohomology space with adjoint coefficients is zero.  

Examples of non-reduced points in the variety of Lie algebra laws have been constructed by
Richardson in \cite{Ric}. Richardson shows there that for the Lie algebra
${\mathfrak s}{\mathfrak l}_2(\C)$ and for $M$ the standard irreducible
${\mathfrak s}{\mathfrak l}_2(\C)$-module of dimension $k+1$,
the semidirect product Lie algebra
${\mathfrak g}:={\mathfrak s}{\mathfrak l}_2(\C)\ltimes M$ is not
rigid if and only if $k=2,4,6,10$. In fact, in these dimensions, there exists a semisimple Lie
algebra of dimension $6,8,10$ and $14$ with ${\mathfrak s}{\mathfrak l}_2(\C)$ as a subalgebra
such that the quotient module identifies with $M$. On the other hand, Richardson shows that
if $\frac{k}{2}>5$ is an odd integer, the second adjoint Lie algebra cohomology space of
${\mathfrak g}$ is non-trivial.
He concludes that if $\frac{k}{2}>5$ is an odd integer, ${\mathfrak g}$ represents a non-reduced point in the variety of
Lie algebra laws.

Going into some more details, Richardson's rigidity result relies on his stability theorem
(joint work with Stanley Page),
see \cite{PagRic}. It says that in case some relative cohomology space
$E^2({\mathfrak g};{\mathfrak s},{\mathfrak g})$ is zero, the subalgebra
${\mathfrak s}\subset{\mathfrak g}$ of the Lie algebra ${\mathfrak g}$ is stable, i.e. all
Lie algebra laws in some neighborhood of ${\mathfrak g}$ contain an isomorphic subalgebra
which has the same brackets with the quotient module. Richardson derives from this theorem
that the above Lie algebra ${\mathfrak g}$ is rigid for $k>10$. The subalgebra ${\mathfrak g}$
in his case is the simple Lie algebra ${\mathfrak s}{\mathfrak l}_2(\C)$ and from here one
can derive that the relative cohomology $E^2({\mathfrak g};{\mathfrak s},{\mathfrak g})$
is zero.

In the Leibniz case, a non-Lie Leibniz algebra with semisimple quotient Lie algebra has necessarily
trivial Leibniz cohomology in degrees $n\geq 2$ (see Proposition
\ref{coh_vanishing_semi_simple}), as follows from Pirashvili's work
\cite{Pir}, see also \cite{FelWag}. Our idea is therefore to take a (non-Lie) Leibniz algebra whose quotient
Lie algebra is of Richardson's type, i.e. a semidirect product of a simple Lie algebra with
an irreducible module. We can show that this kind of Leibniz algebra still has in its adjoint
Leibniz
cohomology the adjoint Lie algebra cohomology of the quotient Lie algebra as a direct factor, see Proposition
\ref{cohomology_double_semi_direct}. This will then imply that in Richardson's setting, the second
Leibniz cohomology of our Leibniz algebra is non-zero.

On the other hand, it is rather straightforward to generalize the proof of the Stability
Theorem to Leibniz algebras as its relies on three applications of the Inverse Function Theorem (see Theorem
\ref{inverse_function_theorem}) using standard
material like Massey products and the coboundary operator. Our stability theorem is Theorem
\ref{stability_theorem}. We show that our example satisfies the hypotheses of this
theorem, namely that the corresponding relative cohomology space
$E^2({\mathfrak h};{\mathfrak s},{\mathfrak h})$ is zero. This is done in Propositions \ref{vanishing} and \ref{computations}.
Let $M_k$ and $I_l$ be finite-dimensional non-trivial irreducible standard left
${\mathfrak s}{\mathfrak l}_2(\C)$-modules of dimensions $k+1$ resp. $l+1$.
Our main theorem reads:

\vspace{.5cm}
\noindent{\bf Theorem.}
{\it The Leibniz algebra ${\mathfrak h}:=I_l\dot{+}(M_k\rtimes{\mathfrak s}{\mathfrak l}_2(\C))$
for two standard irreducible left ${\mathfrak s}{\mathfrak l}_2(\C)$-modules $M_k$ and $I_l$
of highest weights $k=2n$ and $l$ respectively with odd integer $n>5$ and odd $l>2$ is rigid and satisfies
$HL^2({\mathfrak h},{\mathfrak h})\not= 0$. }
\vspace{.5cm}

The structure of the present article is as follows: In a first section, we gather preliminaries
on Leibniz algebras, their modules and semidirect products. The second section is about cohomology
computations. The third section is then devoted to the stability theorem.
The last section concludes the construction of a rigid Leibniz algebra with non-trivial
$HL^2$. In Appendix A, we continue the investigation, started in Section 2.2, of the question for which Lie algebras ${\mathfrak g}$
we have that $H^2({\mathfrak g},{\mathfrak g})=0$ implies $HL^2({\mathfrak g},{\mathfrak g})=0$.\\

\noindent{\bf Acknowledgements:} Bakhrom Omirov thanks Laboratoire de Math\'ematiques Jean Leray for
hospitality during his stay in Nantes in May 2015. Both authors thank J\"org Feldvoss for proofreading the article. 


\section{Preliminaries}   \label{section_preliminaries}


We will always work over a field $k$ of characteristic zero. In geometric situations, we suppose $k$ furthermore to be algebraically closed. Basic material on (right) Leibniz algebras
and their bimodules can be found in \cite{LodPir1}. For left Leibniz algebras, see e.g. \cite{Fel},\cite{FelWag}.

\begin{defi}
A (left) Leibniz algebra is a vector space ${\mathfrak h}$ equipped with a bilinear bracket
$[,]:{\mathfrak h}\times{\mathfrak h}\to{\mathfrak h}$ such that for all $x,y,z\in{\mathfrak h}$
$$[x,[y,z]]\,=\,[[x,y],z]+[y,[x,z]].$$
A morphism of Leibniz algebras is a linear map which respects the brackets.
\end{defi}

Obviously, Lie algebras are examples of Leibniz algebras, but there exist non-Lie Leibniz algebras
(see e.g. Example D in \cite{FelWag}).
We will distinguish Leibniz algebras and Lie algebras in notation by
using ${\mathfrak h}$ for a representative of the former class and ${\mathfrak g}$ for
a representative of the latter class.

Let us recall the notion of a semidirect product of a Lie algebra ${\mathfrak g}$ and a
right ${\mathfrak g}$-module $M$.

\begin{defi}
The semidirect product $M\rtimes{\mathfrak g}$ is the Lie algebra defined on the vector space
$M\oplus{\mathfrak g}$ by
$$[(m_1,x_1),(m_2,x_2)]\,=\,(x_1\cdot m_2-x_2\cdot m_1,[x_1,x_2]).$$
\end{defi}

The action of Lie- or Leibniz algebras (bi-)modules will be written $x\cdot m$ or sometimes
with bracket notation $[x,m]$, following Loday.

\begin{defi}
Let ${\mathfrak h}$ be a (left) Leibniz algebra. A vector space $M$ is called a Leibniz
${\mathfrak h}$-bimodule in case there exist bilinear maps $[,]:{\mathfrak h}\times M\to M$
and $[,]:M\times{\mathfrak h}\to M$ such that for all $x,y\in{\mathfrak h}$ and all $m\in M$
\begin{enumerate}
\item[] (LLM) \quad $[x,[y,m]]\,=\,[[x,y],m]+[y,[x,m]]$,
\item[] (LML) \quad $[x,[m,y]]\,=\,[[x,m],y]+[m,[x,y]]$,
\item[] (MLL) \quad $[m,[x,y]]\,=\,[[m,x],y]+[x,[m,y]]$.
\end{enumerate}
A morphism of Leibniz ${\mathfrak h}$-bimodules is a linear map which respects the
two bracket actions.
\end{defi}

The above three conditions turn up naturally by writing what it means for an abelian
extension to be a Leibniz algebra, see \cite{LodPir1}. For a Lie algebra ${\mathfrak g}$, a left
${\mathfrak g}$-module can be seen as a Leibniz ${\mathfrak g}$-bimodule in two different
ways, namely as a symmetric and as an antisymmetric Leibniz bimodule:

\begin{defi}
\begin{enumerate}
\item[(a)] A Leibniz ${\mathfrak h}$-bimodule $M$ is called symmetric in case for all $x\in{\mathfrak h}$ and all $m\in M$
$$[m,x]\,=\,-[x,m].$$
 \item[(b)] A Leibniz ${\mathfrak h}$-bimodule $M$ is called antisymmetric in case for all
$x\in{\mathfrak h}$ and all $m\in M$
$$[m,x]\,=\,0.$$
\end{enumerate}
\end{defi}

The most important example of an antisymmetric Leibniz ${\mathfrak h}$-bimodule is the
ideal of squares ${\rm Leib}({\mathfrak h})$, also denoted ${\mathfrak h}^{\rm ann}$,
i.e. the ideal of ${\mathfrak h}$ generated by the elements of the
form $[x,x]$ for $x\in{\mathfrak h}$. Indeed, ${\rm Leib}({\mathfrak h})$ becomes a Leibniz
bimodule with respect to the adjoint action and we have
$$[[x,x],y]\,=\,[x,[x,y]]-[x,[x,y]]\,=\,0.$$
The quotient of ${\mathfrak h}$ by ${\rm Leib}({\mathfrak h})$, denoted ${\mathfrak h}_{\rm Lie}$,
is a Lie algebra, and we have an exact sequence of Leibniz algebras
$$0\to {\mathfrak h}^{\rm ann}\to{\mathfrak h}\to{\mathfrak h}_{\rm Lie}\to 0,$$
which is also an abelian extension of Leibniz algebras.

Another important example of an antisymmetric Leibniz ${\mathfrak h}$-module is the
{\it left center} $Z_{\rm left}({\mathfrak h})$ of the Leibniz algebra ${\mathfrak h}$.
It consists by definition of the elements $z\in{\mathfrak h}$ such that for all
$x\in{\mathfrak h}$ we have $[z,x]=0$. By the above, we see that $Z_{\rm left}({\mathfrak h})$
contains the ideal of squares, and the quotient
${\mathfrak h}\,/\,Z_{\rm left}({\mathfrak h})$ is thus a Lie algebra as well.

By quotienting a Leibniz ${\mathfrak h}$-bimodule $M$ by the subbimodule $M^a$ generated by the
elements $[x,m]+[m,x]$ for all $x\in{\mathfrak h}$ and all $m\in M$, one obtains a
symmetric Leibniz bimodule $M^s$, cf \cite{Fel}. The kernel $M^a$ of the projection map
$M\to M^s$ is an antisymmetric bimodule, cf \cite{Fel}. Therefore,
for each Leibniz bimodule $M$, there is a short exact sequence of
Leibniz ${\mathfrak h}$-modules
\begin{equation}  \label{anti_symm_module}
0\to M^a\to M\to M^s\to 0.
\end{equation}

There is also a notion of semidirect product associated to a Lie algebra ${\mathfrak g}$
and a ${\mathfrak g}$-module $I$ which gives a non-Lie Leibniz algebra in case the action is
non-trivial. This notion is due to Kinyon and Weinstein \cite{KinWei}.

\begin{defi}
The hemisemidirect product $I\dot{+}{\mathfrak g}$ is the Leibniz algebra defined on
$I\oplus{\mathfrak g}$ by
$$[(m_1,x_1),(m_2,x_2)]\,=\,(x_1\cdot m_2,[x_1,x_2]).$$
\end{defi}

Note that the hemisemidirect product is split, as ${\mathfrak g}$ embeds as a Leibniz subalgebra in $I\dot{+}{\mathfrak g}$. It constitutes the zero extension in the abelian group of extensions of the Leibniz algebra ${\mathfrak g}$ by the antisymmetric ${\mathfrak g}$-bimodule $I$.

Recall that the Lie algebra ${\mathfrak s}{\mathfrak l}_2(\C)$ admits an irreducible left module $M_k$
of highest weight $k$ for every integer $k\in\N$, $k\geq 0$. The module $M_k$ has dimension $k+1$.
The low dimensional modules are $M_0=\C$ the trivial module, $M_1=\C^2$ the natural module and $M_2={\mathfrak s}{\mathfrak l}_2(\C)$ the adjoint module.
Recall furthermore the Clebsch-Gordan formula (see e.g. \cite{Hum} p. 126):

\begin{prop}  \label{clebsch_gordan}
The tensor product of irreducible left ${\mathfrak s}{\mathfrak l}_2(\C)$-modules $M_m$ and $M_n$ for $m\geq n\in\N$ decomposes into irreducible modules as follows:
$$M_m\otimes M_n\cong M_{n+m}\oplus M_{n+m-2}\oplus\ldots\oplus M_{m-n}.$$
\end{prop}


\section{Leibniz cohomology computations}  \label{section_cohomology}


This section is dedicated to the cohomology of Leibniz algebras, building on the computational methods of \cite{FelWag}, a large part of which is due to Pirashvili \cite{Pir}. 
Later in this section, we will introduce our main example of a Leibniz algebra ${\mathfrak h}$. We will show in this article that ${\mathfrak h}$ is geometrically rigid, but satisfies $HL^2({\mathfrak h},{\mathfrak h})\not= 0$. We introduce in this section the main cohomological tools to assure that $HL^2({\mathfrak h},{\mathfrak h})\not= 0$, while the geometrical rigidity of ${\mathfrak h}$ is the subject of later sections.

\subsection{Leibniz cohomology}

Let ${\mathfrak g}$ be a Lie algebra over a field $k$,
and $M$ a left ${\mathfrak g}$-module. Cohomology theory
associates to ${\mathfrak g}$ two complexes, namely the Chevalley-Eilenberg complex
$$C^*({\mathfrak g},M)\,:=\,(\Hom(\Lambda^*{\mathfrak g},M),d'),$$
and the Leibniz- or Loday complex
$$CL^*({\mathfrak g},M)\,:=\,(\Hom(\bigotimes^*{\mathfrak g},M),d),$$
for the Leibniz cohomology of ${\mathfrak g}$ with values in the symmetric (or antisymmetric)
Leibniz ${\mathfrak g}$-bimodule $M$. The coboundary operator on the complex $C^*({\mathfrak g},M)$ is the standard Chevalley-Eilenberg coboundary operator, see e.g. \cite{HS}. The coboundary operator on $CL^*({\mathfrak g},M)$ is the Leibniz- or Loday coboundary operator
$d:CL^n({\mathfrak g},M)\to CL^{n+1}({\mathfrak g},M)$ defined by
\begin{eqnarray}  \label{coboundary_operator}
(df)(x_1,\dots,x_{n+1}) & := & \sum_{i=1}^n(-1)^{i+1}x_i\cdot f(x_1,\dots,\hat{x}_i,
\dots,x_{n+1})         \\
& + & (-1)^{n+1}f(x_1,\dots,x_n)\cdot x_{n+1}\nonumber    \\
& + & \sum_{1\le i<j\le n+1}(-1)^if(x_1,\dots,\hat{x}_i,\dots,x_ix_j,\dots,x_{n+1}) \nonumber
\end{eqnarray}
for any $f\in CL^n({\mathfrak g},M)$ and all elements $x_1,\dots,x_{n+1}\in{\mathfrak g}$. Leibniz cohomology is more generally defined for any Leibniz algebra ${\mathfrak h}$ and any Leibniz
${\mathfrak h}$-bimodule $M$ with the same coboundary operator $d$, see e.g. \cite{FelWag} for further details.

With values in the symmetric ${\mathfrak g}$-bimodule $M$, the natural epimorphism $\bigotimes^*{\mathfrak g}\to\Lambda^*{\mathfrak g}$ induces a
monomorphism of complexes
$$\varphi:C^*({\mathfrak g},M)\hookrightarrow CL^*({\mathfrak g},M),$$
which is an isomorphism in degree $0$ and $1$. There is then a short exact sequence of complexes
inducing a long exact sequence in cohomology which mediates between Chevalley-Eilenberg and Leibniz cohomology.

It is shown in Lemma 1.5 in \cite{FelWag} that Leibniz cohomology of a Leibniz algebra ${\mathfrak h}$ with values in an antisymmetric ${\mathfrak h}$-bimodule $M^a$ reduces to lower degree cohomology with values in a symmetric ${\mathfrak h}$-bimodule:
\begin{equation}  \label{cohomology_antisymmetric_bimodules}
HL^p({\mathfrak h},M^a)\cong HL^{p-1}({\mathfrak h},\Hom({\mathfrak h},M)^s)
\end{equation}
for all $p\geq 1$.

As usual in cohomology, for a short exact sequence of Leibniz ${\mathfrak h}$-bimodules
$$0\to M'\to M\to M''\to 0,$$
there is a long exact sequence in cohomology
$$\ldots\to HL^n({\mathfrak h},M')\to HL^n({\mathfrak h},M)\to
HL^n({\mathfrak h},M'')\to HL^{n+1}({\mathfrak h},M')\to\ldots,$$
for all $n\geq 0$ and starting with a monomorphism
$$HL^0({\mathfrak h},M')\to HL^0({\mathfrak h},M).$$

\subsection{Leibniz cohomology of rigid Lie algebras}
In this subsection, fix the base field $k$ to be the field $\C$ of complex numbers.

Observe that in general $H^n({\mathfrak g},M)=0$ for one $n$ does not necessarily imply that
$HL^n({\mathfrak g},M)=0$. For example, the trivial Lie algebra ${\mathfrak g}=k$
has $H^2(k,k)=0$ (and is thus Lie-rigid!), but $HL^2(k,k)=k\not= 0$, cf the remark
after Proposition 1 in \cite{FMM}. For another counter-example, a short computation with the Hochschild-Serre spectral sequence \cite{HS}
shows that we have also $H^2({\mathfrak g},{\mathfrak g})=0$ for the direct sum ${\mathfrak g}={\mathfrak s}{\mathfrak l}_2(\C)\oplus\C$, while $HL^2({\mathfrak g},{\mathfrak g})\cong\C$ follows here from the proof of Cor. 3 in \cite{FMM}.

We will consider in this subsection the question for which finite-di\-men\-sio\-nal Lie algebras
${\mathfrak g}\not= k$ the hypothesis
$H^2({\mathfrak g},{\mathfrak g})=0$ implies that $HL^2({\mathfrak g},{\mathfrak g})=0$.
This assertion is true for nilpotent Lie algebras of dimension $\geq 2$, because Th\'eor\`eme 2 of \cite{Dix}
shows that for a non-trivial nilpotent Lie algebra $\dim\,H^2({\mathfrak g},{\mathfrak g})\geq 2$
as its center is non-trivial. The assertion is also true for semisimple Lie algebras
by Proposition \ref{coh_vanishing_semi_simple} below.
The assertion is also true for (non-nilpotent) solvable Lie algebras.
In order to show this, let us elaborate a little on an article of Carles \cite{Car}.

Carles \cite{Car} investigates Lie algebras ${\mathfrak g}$ possessing a codimension one ideal. For these, he shows (in Prop. 2.19) that the dimension $\dim\,\triangle({\mathfrak g})$ of the Lie algebra of derivations $\triangle({\mathfrak g})$ of ${\mathfrak g}$ is greater or equal to the dimension $\dim\,{\mathfrak g}$ of ${\mathfrak g}$. As a corollary (Cor.2.20), this is true for Lie algebras ${\mathfrak g}$ with $[{\mathfrak g},{\mathfrak g}]\not={\mathfrak g}$, because in this case ${\mathfrak g}$ admits a an ideal of codimension one. In a later section, Carles investigates Lie algebras ${\mathfrak g}$ which satisfy $\dim\,\triangle({\mathfrak g})=\dim\,{\mathfrak g}$. He shows (in Prop. 3.1) that any such Lie algebra is algebraic and admits therefore a decomposition ${\mathfrak s}\oplus{\mathfrak u}\oplus{\mathfrak n}$ where ${\mathfrak s}$ is a Levi subalgebra, ${\mathfrak n}$ is the greatest nilpotent ideal and ${\mathfrak u}$ is a subalgebra consisting of ${\rm ad}$-semisimple elements with $[{\mathfrak s}+{\mathfrak u},{\mathfrak u}]=0$
and such that the exterior torus ${\rm ad}({\mathfrak u})$ is algebraic (see Prop. 1.5). In addition, Carles shows in Prop. 3.1 that if furthermore codim $[{\mathfrak g},{\mathfrak g}]>1$, then ${\mathfrak g}$ is complete, i.e. $H^0({\mathfrak g},{\mathfrak g})=H^1({\mathfrak g},{\mathfrak g})=\{0\}$.
In Lemma 5.1, Carles shows that if  $[{\mathfrak g},{\mathfrak g}]\not={\mathfrak g}$,  $\dim\,\triangle({\mathfrak g})\leq \dim\,{\mathfrak g}+\dim\,H^2({\mathfrak g},{\mathfrak g})$.

Let us draw conclusions from these results with respect to the above question.
Let us suppose that ${\mathfrak g}$ is a finite-dimensional solvable Lie algebra with $H^2({\mathfrak g},{\mathfrak g})=0$ and codim $[{\mathfrak g},{\mathfrak g}]>1$. By Carles article, we have on the one hand $\dim\,\triangle({\mathfrak g})\geq \dim\,{\mathfrak g}$ (by Cor. 2.20, because ${\mathfrak g}$ sovable implies $[{\mathfrak g},{\mathfrak g}]\not={\mathfrak g}$ and thus ${\mathfrak g}$ admits an ideal of codimension one). On the other hand, we have $\dim\,\triangle({\mathfrak g})\leq \dim\,{\mathfrak g}$ (by Lemma 5.1, because $H^2({\mathfrak g},{\mathfrak g})=0$). Thus we conclude
$\dim\,\triangle({\mathfrak g})=\dim\,{\mathfrak g}$, which implies then by Prop. 3.1 of Carles that ${\mathfrak g}$ is algebraic and, thanks to codim $[{\mathfrak g},{\mathfrak g}]>1$, ${\mathfrak g}$ is complete, i.e. $H^0({\mathfrak g},{\mathfrak g})=H^1({\mathfrak g},{\mathfrak g})=\{0\}$. This implies that $Z({\mathfrak g})=\{0\}$ and therefore by Theorem \ref{theorem_FMM} below that
$HL^2({\mathfrak g},{\mathfrak g})=\{0\}$.

Let us cite a part of Theorem 2 from \cite{FMM}:

\begin{theo}  \label{theorem_FMM}
Let ${\mathfrak g}$ be a finite-dimensional complex Lie algebra. Then
$H^2({\mathfrak g},{\mathfrak g})$ is a direct factor of
$HL^2({\mathfrak g},{\mathfrak g})$. Furthermore, the supplementary subspace
vanishes in case the center $Z({\mathfrak g})$ is zero.
\end{theo}

\begin{proof}
A proof of this result (as well as further discussion and extensions) is available in \cite{FelWag}.
\end{proof}

Observe that this implies in particular that Richardson's example $M\rtimes{\mathfrak g}$ (see \cite{Ric}) has
non-trivial $HL^2(M\rtimes{\mathfrak g},M\rtimes{\mathfrak g})$.

\begin{cor}  \label{proposition_bakhrom}
A finite dimensional solvable non-nilpotent Lie algebra ${\mathfrak g}$ with $H^2({\mathfrak g},{\mathfrak g})=0$ is the semidirect product of its nilradical ${\mathfrak n}$ and an exterior torus of derivations ${\mathfrak q}$, i.e. ${\mathfrak g}={\mathfrak n}\rtimes{\mathfrak q}$. Furthermore, if $\dim\,{\mathfrak q}>1$,  then $HL^2({\mathfrak g},{\mathfrak g})=\{0\}$.
\end{cor}

\begin{proof}
By the above discussion of the results of Carles in \cite{Car}, it follows that any solvable Lie algebra ${\mathfrak g}$ with $H^2({\mathfrak g},{\mathfrak g})=0$ is
algebraic, i.e., it is isomorphic to the Lie algebra of an algebraic group. As
the algebraicity implies (by Prop. 1.5 in \cite{Car}) the decomposability of the algebra, it
follows that for solvable Lie algebras ${\mathfrak g}$ with $H^2({\mathfrak g},{\mathfrak g})=0$, we have a decomposition
${\mathfrak g}={\mathfrak n}\rtimes{\mathfrak q}$, where ${\mathfrak n}$ is the nilradical of
${\mathfrak g}$ and ${\mathfrak q}$ is an exterior torus of derivations
in the sense of Malcev; that is, ${\mathfrak q}$ is an abelian subalgebra of ${\mathfrak g}$ such that
${\rm ad}(x)$ is semisimple for all $x\in{\mathfrak q}$.

If $\dim\,{\mathfrak q}>1$, it follows from the above discussion before Theorem \ref{theorem_FMM}that $HL^2({\mathfrak g},{\mathfrak g})=\{0\}$.
\end{proof}

Further results in this direction are available in Appendix A.

Ley us now record the Leibniz cohomology of semisimple Lie algebras.
In the context of the following proposition, we assume  that ${\mathfrak g}$ is a semisimple Lie algebra over a field $k$ of characteristic zero.

\begin{prop}  \label{coh_vanishing_semi_simple}
Let $M$ be a finite dimensional left ${\mathfrak g}$-module and $A$ be a
finite dimensional Leibniz ${\mathfrak g}$-bimodule.
Then
$$HL^n({\mathfrak g},M^s)\,=\,0\,\,\,{\rm for}\,\,\,n>0,\,\,\,{\rm and}\,\,\,
HL^n({\mathfrak g},A)\,=\,0\,\,\,{\rm for}\,\,\,n>1.$$
\end{prop}

\begin{prf}
This is the Theorem of Ntolo and Pirashvili \cite{Nto} \cite{Pir}. 
It has been generalized to semisimple Leibniz algebras ${\mathfrak g}$, cf Theorem 4.2 in \cite{FelWag}. 
\end{prf}
 
\subsection{Non-triviality of Leibniz cohomology}

Let $f:{\mathfrak h}\to{\mathfrak q}$ be the quotient morphism
which sends a Leibniz algebras ${\mathfrak h}$ onto its quotient by some two-sided ideal
${\mathfrak k}$, and let $M$ be a Leibniz ${\mathfrak q}$-bimodule. Then $M$ is also a Leibniz
${\mathfrak h}$-bimodule via $f$.
There is a monomorphism of cochain complexes
$$f^*:CL^*({\mathfrak q},M)\to CL^*({\mathfrak h},M),$$
and a quotient complex, called the {\it relative complex}
$$CL^*({\mathfrak h};{\mathfrak q},M)[1]\,:=\,
\coker(f^*:CL^*({\mathfrak q},M)\to CL^*({\mathfrak h},M)).$$
The cohomology spaces are denoted accordingly $HL^*({\mathfrak h};{\mathfrak q},M)$.
The corresponding short exact sequence of complexes induces a long exact sequence in cohomology:

\begin{prop}  \label{long_exact_sequence}
The map $f$ induces a long exact sequence
$$0\to HL^1({\mathfrak q},M)\to HL^1({\mathfrak h},M)\to
HL^1({\mathfrak h};{\mathfrak q},M)\to
HL^2({\mathfrak q},M)\to\ldots $$
\end{prop}

\begin{prf}
This is Proposition 3.1 in \cite{FelWag}.
\end{prf}

We now come to the main example of a Leibniz algebra in our article. Namely, let us consider the semidirect
product Lie algebra $\widehat{\mathfrak g}=M\rtimes{\mathfrak g}$ where ${\mathfrak g}$
is a semisimple Lie algebra (over $\C$) and $M$ is a non-trivial finite-dimensional
irreducible ${\mathfrak g}$-module.
From $\widehat{\mathfrak g}$, we construct as our main object the Leibniz algebra ${\mathfrak h}$ which is
the hemisemidirect product ${\mathfrak h}:=I\dot{+}\widehat{\mathfrak g}$ of
$\widehat{\mathfrak g}$ with ideal of squares $I$ which is another non-trivial finite-dimensional irreducible ${\mathfrak g}$-module. We will need now the above observation that the hemisemidirect product is a split extension.

\begin{prop}  \label{cohomology_double_semi_direct}
Let ${\mathfrak g}$ be a semisimple Lie algebra with non-trivial finite-dimensional
irreducible ${\mathfrak g}$-modules $M$ and $I$.
Then the Leibniz algebra
${\mathfrak h}=I\dot{+}(M\rtimes{\mathfrak g})$ satisfies
$$H^2(\widehat{\mathfrak g},\widehat{\mathfrak g})\hookrightarrow
HL^2({\mathfrak h},{\mathfrak h}),$$
where $\widehat{\mathfrak g}=M\rtimes{\mathfrak g}$.
\end{prop}

\begin{prf}
By Theorem \ref{theorem_FMM}, we have that $H^2(\widehat{\mathfrak g},\widehat{\mathfrak g})$ is a direct factor of $HL^2(\widehat{\mathfrak g},\widehat{\mathfrak g})$.

From now on, as coefficients of the following Leibniz cohomology spaces, we
will consider the short exact sequence of
${\mathfrak h}$-bimodules
$$0\to I\to{\mathfrak h}\to\widehat{\mathfrak g}\to 0.$$
Here the two sided ideal $I$ is an antisymmetric ${\mathfrak h}$-bimodule $I^a$ and
the {\it quotient} ${\mathfrak h}$-bimodule $\widehat{\mathfrak g}$ is in fact a symmetric
${\mathfrak h}$-bimodule where $I$ acts trivially from both sides.    
Therefore, we can apply the preceding constructions. On the other hand, the long
exact sequence for $f:{\mathfrak h}\to\widehat{\mathfrak g}$ of Proposition
\ref{long_exact_sequence} splits, because ${\mathfrak h}$ is the hemisemidirect product of $\widehat{\mathfrak g}$ and $I$. Therefore all connecting homomorphisms are zero and we have a monomorphism for all $n\geq 1$
$$HL^n(\widehat{\mathfrak g},\widehat{\mathfrak g})\hookrightarrow HL^n({\mathfrak h},\widehat{\mathfrak g}).$$

Recall the construction of the connecting homomorphism
$\partial:HL^2({\mathfrak h},\widehat{\mathfrak g})\to HL^3({\mathfrak h},I^a)$
in the long exact sequence for the Leibniz cohomology of ${\mathfrak h}$ with
values in the short exact sequence of coefficients
$$0\to I\to{\mathfrak h}\to\widehat{\mathfrak g}\to 0.$$
We claim that the subspace $HL^2(\widehat{\mathfrak g},\widehat{\mathfrak g})\subset HL^2({\mathfrak h},\widehat{\mathfrak g})$ is in the kernel of $\partial$. This is clear, because lifting a $2$-cocycle
$c\in CL^2(\widehat{\mathfrak g},\widehat{\mathfrak g})$ to a cochain in 
$CL^2({\mathfrak h},{\mathfrak h})$, it will remain a cocycle and thus the preimage of its coboundary in  $CL^3({\mathfrak h},I^a)$ is zero. Therefore we have an epimorphism
$$HL^2({\mathfrak h},{\mathfrak h})\twoheadrightarrow HL^2(\widehat{\mathfrak g},\widehat{\mathfrak g}).$$
This ends the proof of the proposition, because all cohomology spaces are $\C$-vector spaces and thus the above epimorphism splits.
\end{prf}

\subsection{Vanishing of some relative cohomology groups}

In this subsection, we consider a type of relative cohomology which will be useful
in the construction of a rigid Leibniz algebra ${\mathfrak h}$ whose
$HL^2({\mathfrak h},{\mathfrak h})$ is not zero.

Let ${\mathfrak h}$ be a finite-dimensional Leibniz algebra and ${\mathfrak s}$ be
a subalgebra. Let ${\mathfrak h}^n$ be the $n$-fold Cartesian product and denote by
$F_n$ the subset of ${\mathfrak h}^n$ consisting of all $(x_1,\ldots,x_n)$ satisfying
the following condition: There exists at most one index $i$ such that
$x_i\notin{\mathfrak s}$. Let $M$ be a Leibniz ${\mathfrak h}$-module.
We define $F^n({\mathfrak h};{\mathfrak s},M)$ to be the vector space of all
multilinear maps $\phi:F_n\to M$. Concretely for $n=2$,
$\phi\in F^2({\mathfrak h};{\mathfrak s},M)$ means that $\phi$ is defined on
$({\mathfrak h}\times{\mathfrak s})\cup({\mathfrak s}\times{\mathfrak h})\cup
({\mathfrak s}\times{\mathfrak s})$.
We set
$$F({\mathfrak h};{\mathfrak s},M)\,=\,\bigoplus_{n\geq 0}F^n({\mathfrak h};{\mathfrak s},M).$$
The coboundary operator $d$ for Leibniz cohomology is defined as in Equation (\ref{coboundary_operator}). One
checks that it is well-defined on $F({\mathfrak h};{\mathfrak s},M)$. We denote by
$$E({\mathfrak h};{\mathfrak s},M)\,=\,\bigoplus_{n\geq 0}E^n({\mathfrak h};{\mathfrak s},M)$$
the cohomology of this complex. We will need this cohomology with adjoint coefficients,
i.e. with values in the ${\mathfrak h}$-bimodule ${\mathfrak h}$, and consider thus $E({\mathfrak h};{\mathfrak s},{\mathfrak h})$.

Let us first analyze cocycles $\phi\in F^2({\mathfrak h};{\mathfrak s},{\mathfrak h})$.
For this, let $W$ be a supplementary subspace of ${\mathfrak s}$ in ${\mathfrak h}$.
For our main application, we will have $W=M\oplus I$.
The condition
$d\phi=0$ means explicitly for all $x,y,z\in{\mathfrak h}$
$$[x,\phi(y,z)]-[y,\phi(x,z)]-[\phi(x,y),z]-\phi([x,y],z)-\phi(y,[x,z])+\phi(x,[y,z])\,=\,0.$$
A priori, $\phi$ is a map from $({\mathfrak h}\times{\mathfrak s})\cup({\mathfrak s}\times{\mathfrak h})\cup({\mathfrak s}\times{\mathfrak s})$ to ${\mathfrak h}$. Writing ${\mathfrak h}=W\oplus{\mathfrak s}$ as vector spaces,
$\phi$ splits into three components $f_1:{\mathfrak s}\otimes W\to{\mathfrak h}$,
$f_2:W\otimes{\mathfrak s}\to{\mathfrak h}$ and
$f_3:{\mathfrak s}\otimes {\mathfrak s}\to{\mathfrak h}$. In terms of these components, the
cocycle condition reads
\begin{eqnarray}
\forall x,y\in{\mathfrak s},\forall z\in W:\,\,\,[x,f_1(y,z)]-[y,f_1(x,z)]-[f_3(x,y),z]-
\nonumber\\
f_1([x,y],z)-f_1(y,[x,z])+f_1(x,[y,z])=0,  \label{****}\\
\forall x,z\in{\mathfrak s},\forall y\in W:\,\,\,[x,f_2(y,z)]-[y,f_3(x,z)]-[f_1(x,y),z]-
\nonumber\\
f_2([x,y],z)-f_2(y,[x,z])+f_1(x,[y,z])=0,  \label{**}\\
\forall y,z\in{\mathfrak s},\forall x\in W:\,\,\,[x,f_3(y,z)]-[y,f_2(x,z)]-[f_2(x,y),z]-
\nonumber\\
f_2([x,y],z)-f_1(y,[x,z])+f_2(x,[y,z])=0.\label{***}
\end{eqnarray}

Similar conditions hold for two components in $W$ and one in ${\mathfrak s}$, but these
simplify, because we now impose that $[W,W]=0$ and $[{\mathfrak s},W]\subset W$ and
$[W,{\mathfrak s}]\subset W$. Recall furthermore that by definition of
$F({\mathfrak h};{\mathfrak s},M)$, we have $f|_{W\otimes W}=0$. The conditions read then:
\begin{eqnarray}
\forall x,y\in W,\forall z\in {\mathfrak s}:\,\,\,[x,f_2(y,z)]+[y,f_2(x,z)]\,=\,0 \nonumber \\
\forall x,z\in W,\forall y\in {\mathfrak s}:\,\,\,[x,f_1(y,z)]-[f_2(x,y),z]\,=\,0 \label{*}\\
\forall y,z\in W,\forall x\in {\mathfrak s}:\,\,\,-[y,f_1(x,z)]-[f_1(x,y),z]\,=\,0\nonumber
\end{eqnarray}
We will use these identities in the proof of the following proposition in order to
express a cocycle $\phi\in F^2({\mathfrak h};{\mathfrak s},M)$ in terms of cocycles
in ordinary Leibniz cohomology.

Before we come to the proposition, let us also record the coboundary relations. A priori, a cochain
$\psi\in F^1({\mathfrak h};{\mathfrak s},{\mathfrak h})$ is a map from ${\mathfrak s}\cup{\mathfrak h}$ to ${\mathfrak h}$. Writing ${\mathfrak h}=W\oplus{\mathfrak s}$, $\psi$ splits into two components
$g_1:{\mathfrak s}\to{\mathfrak h}$ and $g_2:W\to{\mathfrak h}$. The coboundary condition
$d\psi(x,y)=[x,\psi(y)]+[\psi(x),y]-\psi([x,y])$ splits into the coboundary condition for $g_1$ (for $x,y\in{\mathfrak s}$) and the conditions
\begin{eqnarray*}
\forall x\in W,\forall y\in{\mathfrak s}:\,\,\,[x,g_1(y)]+[g_2(x),y]-g_2([x,y]) \\
\forall x\in {\mathfrak s},\forall y\in W:\,\,\,[x,g_2(y)]+[g_1(x),y]-g_2([x,y]),
\end{eqnarray*}
using as above that $[W,{\mathfrak s}]\subset W$ and $[W,{\mathfrak s}]\subset W$.
Observe that the $1$-coboundary equation above has only two terms in case $x,y$ are not both in ${\mathfrak s}$ and $g_1=0$. Therefore the coboundary condition for $g_2$ resembles in this case a $0$-coboundary condition with values in a $\Hom$-space.
We will use this observation in the proof of the following proposition.

\begin{prop}   \label{vanishing}
Let ${\mathfrak h}$ be a Leibniz algebra with subalgebra ${\mathfrak s}$.
Let $M$ and $I$ be subspaces of ${\mathfrak h}$ such that
${\mathfrak h}={\mathfrak s}\oplus M\oplus I$
as a vector space. We suppose that $M\oplus I$ is an ideal of ${\mathfrak h}$,
and in this way, we will consider $M$ and $I$ as ${\mathfrak s}$-modules. Suppose furthermore
\begin{enumerate}
\item[(a)] $[M,M]=0$, $[{\mathfrak s},M]\subset M$, $[M,{\mathfrak s}]\subset M$,
and the bracket of ${\mathfrak s}$ and $M$ is symmetric:
$$\forall s\in {\mathfrak s},\forall w\in M:\,\,\,[s,w]\,=\,-[w,s],$$
furthermore $[I,{\mathfrak s}]=0$, $[{\mathfrak s},I]\subset I$ and
$I=Z_{\rm left}({\mathfrak h})={\rm Leib}({\mathfrak h})$, and suppose that the ${\mathfrak s}$-annihilators of the ${\mathfrak s}$-modules $M$ and $I$ are trivial,
\item[(b)] $HL^2({\mathfrak s},{\mathfrak h})=0$,
\item[(c)] $HL^1({\mathfrak h},\Hom(M^s\oplus I^s,{\mathfrak h}^s)=0$,
\item[(d)] $HL^1({\mathfrak s},\Hom(M,I)^a)=0$,
\item[(e)] $HL^1({\mathfrak s},\Hom(I^{\rm triv},M\oplus I))=0$, where ${\mathfrak s}$ acts trivially on $I^{\rm triv}$.
\end{enumerate}
Then we have $E^2({\mathfrak h};{\mathfrak s},{\mathfrak h})=0$.
\end{prop}

\begin{prf}
Let $\phi\in F^2({\mathfrak h};{\mathfrak s},{\mathfrak h})$ with $d\phi=0$.
First of all, it follows from condition {\it (a)} that the three components $f_1$, $f_2$
and $f_3$ of $\phi$ satisfy the above mentioned equations.
Observe that as $\phi:{\mathfrak h}^{\otimes 2}\to{\mathfrak h}$ is a cocycle, its
restriction $f_3$ to ${\mathfrak s}\times{\mathfrak s}$ is a $2$-cocycle on ${\mathfrak s}$
with values in the Leibniz ${\mathfrak s}$-bimodule ${\mathfrak h}$. It follows then
from {\it (b)} that there exists $\psi\in F^1({\mathfrak h};{\mathfrak s},{\mathfrak h})$
such that $\phi_1=\phi-d\psi$ vanishes on ${\mathfrak s}\times{\mathfrak s}$.
This means that for the three components of $\phi_1$
(for which we will not introduce new notations!), we may erase $f_3$ from the equations.

Equation \ref{****} for $f_1:{\mathfrak s}\otimes(M\oplus I)\to{\mathfrak h}$ reads then for all $x,y\in{\mathfrak s}$ and all $z\in M\oplus I$:
$$[x,f_1(y,z)]-[y,f_1(x,z)]-f_1([x,y],z)-f_1(y,[x,z])+f_1(x,[y,z])\,=\,0.$$
This can be interpreted as the cocycle identity for $1$-cocycles on ${\mathfrak s}$ with values in $\Hom(M^s\oplus I^s,{\mathfrak h}^s)$. Indeed, observe that in order to interprete it as a cocycle identity, one has to swich the left action by $y$ into a right action which imposes to view $I$, $M$ and ${\mathfrak h}$ as symmetric ${\mathfrak s}$-bimodules.
By condition {\it (c)}, there exists $\psi_1\in F^1({\mathfrak h};{\mathfrak s},{\mathfrak h})$
such that $\phi_2=\phi_1-d\psi_1$ vanishes
on ${\mathfrak s}\times{\mathfrak s}$ and on ${\mathfrak s}\times (M\oplus I)$.
Observe that the coboundary condition for the $1$-cocycle with values in the $\Hom$-space is indeed a coboundary condition for the $2$-cocycle as the latter contains only two terms.
In conclusion, this means that for the three components of $\phi_1$,
we may erase $f_1$ and $f_3$ from the equations.

Now consider $f_2:(M\oplus I)\otimes{\mathfrak s}\to{\mathfrak h}$ (where we again omit to introduce new notations and thus $f_2$ is the corresponding component of $\phi_2$). Observe that Equation (\ref{*}) shows that
$\im(f_2)\subset M\oplus I$, because there is no non-trivial element of ${\mathfrak s}$ which
annihilates all elements of $M\oplus I$ by condition {\it (a)}.
Equation (\ref{**}) reads now for all $x,z\in{\mathfrak s}$ and all $y\in M\oplus I$
\begin{equation}  \label{+}
[x,f_2(y,z)]-f_2([x,y],z)-f_2(y,[x,z])=0.
\end{equation}
Equation (\ref{***}) reads then for all $y,z\in{\mathfrak s}$ and all $x\in M\oplus I$
$$-[y,f_2(x,z)]-[f_2(x,y),z]-f_2([x,y],z)+f_2(x,[y,z])=0.$$
Let us exchange $x$ and $y$ such that the equation becomes for all $x,z\in{\mathfrak s}$ and all $y\in M\oplus I$
\begin{equation}  \label{++}
-[x,f_2(y,z)]-[f_2(y,x),z]-f_2([y,x],z)+f_2(y,[x,z])=0.
\end{equation}
We now restrict to $y\in M$ and use that $M$ is a symmetric ${\mathfrak s}$-bimodule in order
to change $f_2([y,x],z)$ into $-f_2([x,y],z)$. As a result, the sum of Equations (\ref{+}) and (\ref{++})
gives for all $x,z\in{\mathfrak s}$ and all $y\in M$
$$[f_2(y,x),z]=0.$$
This means that this component $f_2:M\otimes{\mathfrak s}\to M\oplus I$ of $f_2$ has values in the left center which we supposed to be $I$ by condition {\it (a)}.
Now Equation (\ref{+}) is exactly the cocycle identity for the $1$-cochain
$f_2$ with values in the antisymmetric ${\mathfrak s}$-bimodule $\Hom(M,I)^a$.
By condition {\it (d)}, there exists $\psi_2\in F^1({\mathfrak h};{\mathfrak s},{\mathfrak h})$
such that $\phi_3=\phi_2-d\psi_2$ vanishes on ${\mathfrak s}\times{\mathfrak s}$, on ${\mathfrak s}\times (M\oplus I)$ and on $M\times{\mathfrak s}$.

The last step is then to investigate the component of $f_2$ which is $f_2:I\otimes{\mathfrak s}\to M\oplus I$. As $I$ is an antisymmetric ${\mathfrak s}$-bimodule by condition {\it (a)}, Equation
(\ref{++}) simplifies to
$$-[x,f_2(y,z)]-[f_2(y,x),z]+f_2(y,[x,z])=0$$
for all $x,z\in{\mathfrak s}$ and all $y\in I$. This is the $1$-cocycle identity for ${\mathfrak s}$ with values in the ${\mathfrak s}$-bimodule $\Hom(I^{\rm triv},M\oplus I)$, where the elements of ${\mathfrak s}$ do not act on the domain $I$, but only on the codomain $M\oplus I$.
By condition {\it (e)}, there exists $\psi_3\in F^1({\mathfrak h};{\mathfrak s},{\mathfrak h})$
such that $\phi_4=\phi_3-d\psi_3$ vanishes on ${\mathfrak s}\times{\mathfrak s}$, on ${\mathfrak s}\times (M\oplus I)$ and on $(M\oplus I)\times{\mathfrak s}$, i.e. is identically zero.
Therefore $\phi=d\psi+d\psi_1+d\psi_2+d\psi_3$ as was to be shown.
\end{prf}

Let us now apply this proposition to our very special Leibniz algebra.
Consider again a semidirect
product Lie algebra $\widehat{\mathfrak g}=M\rtimes{\mathfrak g}$ where ${\mathfrak g}$
is a simple Lie algebra (over $\C$) and $M$ is a non-trivial finite-dimensional
irreducible ${\mathfrak g}$-module.
From $\widehat{\mathfrak g}$, we construct a Leibniz algebra ${\mathfrak h}$ which is
the hemisemidirect product ${\mathfrak h}:=I\dot{+}\widehat{\mathfrak g}$ of
$\widehat{\mathfrak g}$ with ideal of squares $I$ which is another
non-trivial finite-dimensional irreducible ${\mathfrak g}$-module.

As subalgebra ${\mathfrak s}$, we take the simple Lie algebra ${\mathfrak g}$.
The supplementary subspace is thus $M\oplus I$.
Thus the pair $({\mathfrak h},{\mathfrak s})$ satisfies conditions  {\it (a)} to {\it (c)}
of the proposition, using for conditions {\it (b)} and {\it (c)} Proposition \ref{coh_vanishing_semi_simple} as $\Hom(M^s\oplus I^s,{\mathfrak h}^s)$ is a symmetric ${\mathfrak g}$-bimodule. In order to meet conditions {\it (d)} and {\it (e)}, we will have to restrict the irreducible
${\mathfrak g}$-modules $M$ and $I$ of the construction. Namely, we suppose that
${\mathfrak g}={\mathfrak s}{\mathfrak l}_2(\C)$ and that $M$ and $I$ are the standard irreducible
${\mathfrak s}{\mathfrak l}_2(\C)$-modules $M=M_k$ and $I=I_l$ of dimensions $\dim(M_k)=k+1$
and $\dim(I_l)=l+1$ and highest weights $k$ resp. $l$.

\begin{prop}          \label{computations}
Assume that $\dim(M_k)=k+1$ and $\dim(I_l)=l+1>3$ are such that $k$ is even and $l$ is odd.

Then
\begin{enumerate}
\item[(d)] $HL^1({\mathfrak s}{\mathfrak l}_2(\C),\Hom(M_k,I_l)^a)=0$,
\item[(e)] $HL^1({\mathfrak s}{\mathfrak l}_2(\C),\Hom(I^{\rm triv}_l,M_k\oplus I_l))=0$, where ${\mathfrak s}{\mathfrak l}_2(\C)$ acts trivially on $I^{\rm triv}_l$.
\end{enumerate}
\end{prop}

\begin{prf}
Observe that we have isomorphisms $(M_k)^*\cong M_k$
and $(I_l)^*\cong I_l$, because they are the only simple ${\mathfrak s}{\mathfrak l}_2(\C)$-modules in these dimensions.

The main idea is to translate $1$-cohomology with values in an antisymmetric
${\mathfrak s}{\mathfrak l}_2(\C)$-bimodule into $0$-cohomology with values in the homomorphisms, see the isomorphism in Equation (\ref{cohomology_antisymmetric_bimodules}).
We have
$$HL^1({\mathfrak s}{\mathfrak l}_2(\C),\Hom(M_k,I_l)^a)=HL^0({\mathfrak s}{\mathfrak l}_2(\C),\Hom({\mathfrak s}{\mathfrak l}_2(\C),\Hom(M_k,I_l)^s)).$$
This latter space is the space of ${\mathfrak s}{\mathfrak l}_2(\C)$-morphisms from ${\mathfrak s}{\mathfrak l}_2(\C)$ to $\Hom(M_k,I_l)$. As ${\mathfrak s}{\mathfrak l}_2(\C)$ is simple, there can be non-trivial contributions only if $\Hom(M_k,I_l)\cong(M_k)^*\otimes I_l\cong M_k\otimes I_l$
contains a factor isomorphic to
${\mathfrak s}{\mathfrak l}_2(\C)$. The Clebsch-Gordan decomposition formula of the tensor product
in Proposition \ref{clebsch_gordan}
shows that there can be only components of even dimension in the decomposition, because we have chosen $k,l$ such that one of them is even and the other one is odd. Therefore, the space of invariants is zero.

For condition {\it (e)}, we observe that $\Hom(I^{\rm triv}_l,M_k\oplus I_l)=\Hom(I_l^{\rm triv},M_k)\oplus\Hom(I_l^{\rm triv},I_l)$ decomposes into a sum of tensor products $\C\otimes_{\C} M_k$ and $\C\otimes_{\C}I_l$, because the trivial module $I_l^{\rm triv}$ decomposes into a direct sum of trivial modules. Thus the space $HL^1({\mathfrak s}{\mathfrak l}_2(\C),\Hom(I^{\rm triv}_l,M_k\oplus I_l))$ decomposes into a direct sum of spaces of the form $HL^1({\mathfrak s}{\mathfrak l}_2(\C),M_k)$ and $HL^1({\mathfrak s}{\mathfrak l}_2(\C),I_l)$. Furthermore, we have that $HL^1({\mathfrak s}{\mathfrak l}_2(\C),M_k)=0$  by Proposition \ref{coh_vanishing_semi_simple}, because $M_k$ is a symmetric
${\mathfrak s}{\mathfrak l}_2(\C)$-bimodule, and
$$HL^1({\mathfrak s}{\mathfrak l}_2(\C),I_l)=HL^0({\mathfrak s}{\mathfrak l}_2(\C),\Hom({\mathfrak s}{\mathfrak l}_2(\C),I_l)),$$
which is zero, because $l+1>3$.
 \end{prf}


\section{A stability theorem}   \label{section_stability}


In this section, we follow closely \cite{PagRic} in order to show a stability theorem
for Leibniz subalgebras of a given Leibniz algebra. We work over the complex numbers $\C$.

Let $V$ be a finite dimensional complex vector space.
Let ${\mathcal M}$ be the algebraic variety of all Leibniz algebra structures on $V$, defined in
$\Hom(V^{\otimes 2},V)$ by the quadratic equations corresponding to the left Leibniz identity,
i.e. for $\mu\in\Hom(V^{\otimes 2},V)$, we require
$$\mu\circ\mu(x,y,z)\,:=\,\mu(x,\mu(y,z))-\mu(\mu(x,y),z)-\mu(y,\mu(x,z))\,=\,0$$
for all $x,y,z\in V$.

Two Leibniz algebra structures $\mu$ and $\mu'$ give rise to isomorphic Leibniz algebras
$(V,\mu)$ and $(V,\mu')$ in case there exists $g\in\Gl(V)$ such that $g\cdot\mu=\mu'$,
where the action of $\Gl(V)$ on ${\mathcal M}$ is defined by
$$g\cdot\mu(x,y)\,:=\,g(\mu(g^{-1}(x),g^{-1}(y)))$$
for all $x,y\in V$.

Let ${\mathfrak h}=(V,\mu)$ be a Leibniz algebra on $V$, and consider a Leibniz subalgebra
${\mathfrak s}$ of ${\mathfrak h}$. We will denote the complex subspace of $V$
corresponding to ${\mathfrak s}$ by $S$. Let $W$ be a supplementary subspace of $S$ in $V$.
For an element $\phi\in CL^2({\mathfrak h},V)$, we denote
by $r(\phi)$ the restriction to the union of $S\otimes S$, $S\otimes W$ and $W\otimes S$.
Let
$$N_1\,:=\,\{(g,m)\in\Gl(V)\times {\mathcal M}\,|\,r(g\cdot m)=r(m)\},$$
denote by ${\rm proj}_{{\mathcal M}}:\Gl(V)\times {\mathcal M}\to {\mathcal M}$ the projection map and let
$p_1:N_1\to {\mathcal M}$ be the restriction of ${\rm proj}_{{\mathcal M}}$ to $N_1$.

\begin{defi}
The subalgebra ${\mathfrak s}$ is called a stable subalgebra of ${\mathfrak h}$
if $p_1$ maps every neighborhood of $(1,\mu)\in N_1$
onto a neighborhood of $\mu$ in ${\mathcal M}$.
\end{defi}

\begin{rem}
\begin{enumerate}
\item[(a)] This is not the original definition of stable subalgebra in \cite{Ric}, but
the strong version of stability which permits Page and Richardson to show the strengthened
stability theorem at the end of their paper.
\item[(b)] The definition implies that if ${\mathfrak h}_1=(V,\mu_1)$ is a Leibniz algebra
sufficiently near to ${\mathfrak h}$, then ${\mathfrak h}_1$ is isomorphic to a Leibniz
algebra ${\mathfrak h}_2=(V,\mu_2)$ with the following property: for all $s\in{\mathfrak s}$
and all $x\in V$, we have $\mu_2(x,s)=\mu(x,s)$ and $\mu_2(s,x)=\mu(s,x)$.
\end{enumerate}
\end{rem}

The stability theorem below asserts that under certain cohomological conditions a given
 subalgebra ${\mathfrak s}$ is stable. The proof is based on the inverse and implicit
function theorem,
in the following algebraic geometric form. All algebraic varieties are considered
to be complex and equipped with the Zariski topology.
Recall that a point $x$ of an algebraic variety
$X$ is called {\it simple} in case $\dim(X)=\dim(T_xX)$.

\begin{theo}  \label{inverse_function_theorem}
Let $f:X\to Y$ be a morphism of algebraic varieties, let $x\in X$ and $y=f(x)\in Y$.
Suppose that $x$ (resp. $y$) is a simple point of $X$ (resp. $Y$) and that the differential
$d_xf:T_xX\to T_yY$ is surjective.
\begin{enumerate}
\item[(a)] If $U$ is a neighborhood of $x$ in $X$, then $f(U)$ is a neighborhood of $y$
in $Y$.
\item[(b)] There is exactly one irreducible component $X_1$ of $f^{-1}(y)$ which contains $x$.
\end{enumerate}
Moreover, $x$ is a simple point of $X_1$ and $T_xX_1=\ker(d_xf)$.
\end{theo}

\begin{prf} This is Proposition 8.1 in \cite{PagRic}.
\end{prf}

\begin{theo}[Stability Theorem] \label{stability_theorem}
Let ${\mathfrak h}=(V,\mu)$ be a finite-dimensional complex Leibniz algebra, and let
${\mathfrak s}$ be a subalgebra of ${\mathfrak h}$ such that $E^2({\mathfrak h},{\mathfrak s},
V)=0$. Then ${\mathfrak s}$ is a stable subalgebra of ${\mathfrak h}$.
\end{theo}

\begin{prf}
The proof procedes as in Section 11 of \cite{PagRic} by three distinct applications of the inverse
function theorem Theorem \ref{inverse_function_theorem}.

\noindent{\bf (1)} For $\phi\in\Hom(V^{\otimes 2},V)$, let $\tau(\phi)$ denote the
restriction of $\phi$ to $F_2$. Let
$$Z^2\,:=\,\{\theta\in F^2({\mathfrak h},{\mathfrak s},V)\,|\,d\theta=0\}$$
and let $C:=\{\phi\in\Hom(V^{\otimes 2},V)\,|\,\tau(\phi)\in Z^2\}$. Then
$ZL^2({\mathfrak h},V)\subset C$. Denote by $C_0:=d(C)$ and let $C_1$ be a supplementary subspace to
$C_0$ in $B^3(L,V)$ and $C_2$ be a supplementary subspace to $B^3(L,V)$ in $\Hom(V^{\otimes 3},V)$.
Let $\pi_1:\Hom(V^{\otimes 2},V)\to C_1$ be the projection with kernel $C_0\oplus C_2$.
Put
$${\mathcal M}_1:=\{\phi\in\Hom(V^{\otimes 2},V)\,|\,\pi_1(\phi\circ\phi)=0\}.$$
Then ${\mathcal M}\subset {\mathcal M}_1$, because for elements $\mu\in {\mathcal M}$, we have $\mu\circ\mu=0$.

Let $f:\Hom(V^{\otimes 2},V)\to\Hom(V^{\otimes 3},V)$ be defined by $f(\phi)=\phi\circ\phi$.
$f$ is a polynomial mapping.
As in the proof of 6.2 in \cite{PagRic}, we show that $d_{\mu}f=-d$, i.e. that the differential of $f$ at $\mu$ is equal to the Leibniz coboundary operator.
This is true since the Massey product $\circ$ and the Leibniz coboundary operator share the
same relations as for Lie algebras.
Let $F:\Hom(V^{\otimes 2},V)\to C_1$ be the composition $\pi_1\circ f$. Since $\pi_1$ is
linear, $d_{\mu}F=-\pi_1\circ d$. Hence $d_{\mu}F$ is surjective, because $C_1$ consists of
coboundaries ($C_1$ are exactly the coboundaries which do not come from $C$ !), and
$\ker(d_{\mu}F)=C$. By Theorem \ref{inverse_function_theorem}, there is exactly one
irreducible component ${\mathcal M}'$ of ${\mathcal M}_1$ which contains $\mu$. Moreover, $\mu$ is a simple
point of ${\mathcal M}'$ and $T_{\mu}{\mathcal M}'=C$.

\noindent{\bf (2)} Let $E$ be a supplementary subspace of $C$ in $\Hom(V^{\otimes 2},V)$ and let $E'=\tau(E)$.
We have a direct sum decomposition $F^2({\mathfrak h},{\mathfrak s},V)=Z^2\oplus E'$ in
cocycles $Z^2$ and the rest $E'$. Let $\pi_Z:F^2({\mathfrak h},{\mathfrak s},V)\to Z^2$
be the projection with kernel $E'$. Let $\gamma:\Gl(V)\to\Hom(V^{\otimes 2},V)$ be the map
$\gamma(g)=g\cdot\mu$. Taking $g=1+t\psi+{\mathcal O}(t^2)$, one obtains
$$d_1\gamma(\psi)(x,y)=\frac{d}{dt}\left.\,g(\mu(g^{-1}(x),g^{-1}(y)))\right|_{t=0}
=\psi([x,y])-[\psi(x),y]-[x,\psi(y)],$$
for all $x,y\in V$ (with $\mu(x,y)=[x,y]$). Thus $d_1\gamma(\psi)=-d\psi$.

Let $\eta:\Gl(V)\times {\mathcal M}'\to\Hom(V^{\otimes 2},V)$ be
defined by $\eta(g,m)=g\cdot m$ and let $\beta:\Gl(V)\times {\mathcal M}'\to Z^2$ be the composition
$\pi_Z\circ\tau\circ\eta$. Using the above computation for $d_1\gamma$, one obtains for $d_{(1,\mu)}\eta$
$$d_{(1,\mu)}\eta(\psi,\phi)\,=\,-d\psi+\phi.$$
Therefore $d_{(1,\mu)}\beta(\psi,\phi)\,=\,-d\tau(\psi)+\pi_Z(\tau(\phi))$. As $T_{\mu}{\mathcal M}'=C$
and $\tau(C)=Z^2$, it follows that $d_{(1,\mu)}\beta$ is surjective and we can again
apply Theorem \ref{inverse_function_theorem}: There exists exactly one irreducible component
$N'$ of $\beta^{-1}(\tau(\mu))$ containing $(1,\mu)$ and $(1,\mu)$ is a simple point of
$N'$ with $T_{(1,\mu)}N'=\ker(d\beta_{(1,\mu)})$.

\noindent{\bf (3)} Let $\pi:N'\to {\mathcal M}'$ denote the restriction of the projection
$${\rm proj}_{{\mathcal M}'}:\Gl(V)\times {\mathcal M}'\to {\mathcal M}';$$
We claim that if $E^2({\mathfrak h},{\mathfrak s},V)=0$, the differential
$d_{(1,\mu)}\pi:T_{(1,\mu)}N'\to T_{\mu}{\mathcal M}'$ is surjective. Indeed, let $\phi\in T_{\mu}{\mathcal M}'$.
By definition of $C$, $\tau(\phi)\in Z^2=B^2({\mathfrak h},{\mathfrak s},V)$ by assumption.
Thus there exists $\theta\in\Hom(S,V)$ such that $d\theta=\tau(\phi)$. Since $\tau$ is surjective,
there exists $\psi\in\Hom(V,V)$ (extend by zero for example!) such that $\tau(\psi)=\theta$;
therefore $d\tau(\psi)=\tau(\phi)=\pi_Z(\tau(\phi))$. Consequently
$$d_{(1,\mu)}\beta(\psi,\phi)\,=\,-d\tau(\psi)+\pi_Z(\tau(\phi))\,=\,0.$$
It follows that $(\psi,\phi)\in T_{(1,\mu)}N'=\ker(d_{(1,\mu)}\beta)$. Since $d_{(1,\mu)}\pi(\psi,\phi)=\phi$,
this shows that $d_{(1,\mu)}\pi$ is surjective. This will be used later on for the third application of
the inverse function theorem.

\noindent{\bf (4)} We claim now that there exists a neighborhood $U_3$ of $(1,\mu)$ in $\Gl(V)\times {\mathcal M}$ such that
if $(g,m)\in U_3$ and if $\pi_Z(\tau(g\cdot m))=\tau(\mu)$, then $\tau(g\cdot m)=\tau(\mu)$.
Indeed, let $(g,m)\in\Gl(V)\times {\mathcal M}$ with $\pi_Z(\tau(g\cdot m))=\tau(\mu)$. Then
$g\cdot m=\mu+a+b$ with $a\in\ker(\tau)$ and $b\in E$ (because of the direct sum decomposition
$\Hom(V^{\otimes 2},V)=C\oplus E$ which maps under $\tau$ to $Z^2\oplus E'$, and the $Z^2$-component
is fixed to be $\tau(\mu)$). The restriction of $d\circ\tau$ to $E$ is a monomorphism
(by definition of $C$ and $E$). As $g\cdot m$ is a Leibniz law, we have
\begin{eqnarray*} 0&=& (g\cdot m)\circ(g\cdot m)\,=\,(\mu+a+b)\circ(\mu+a+b) \\
&=& -da-db+a\circ a+b\circ b +a\circ b + b\circ a.
\end{eqnarray*}
As $a\in\ker(\tau)$, $\tau(a\circ a)=\tau(da)=0$. Applying $\tau$ to the
preceeding equation, we obtain:
$$0\,=\,-d\tau(b)+\tau((a+b)\circ b)+\tau(b\circ a).$$
Let $\pi_C:\Hom(V^{\otimes 2},V)\to C$ be the projection with kernel $E$. We note that $\pi_C(a+b)=a$.
Thus the above equation can be rewritten as:
$$0\,=\,-d\tau(b)+\tau((a+b)\circ b)+\tau(b\circ\pi_C(a+b)).$$
For any $x\in\Hom(V^{\otimes 2},V)$, let $\lambda_x:E\to\Hom(V^{\otimes 3},V)$ be the linear map
defined by
$$\lambda_x(\phi)\,:=\,-d\tau(\phi)+\tau(x\circ\phi)+\tau(\phi\circ\pi_C(x)).$$
Then $\lambda_0=d\circ\tau$ is a monomorphism. Thus there exists a neighborhood $J$ of $0$
in $\Hom(V^{\otimes 2},V)$ such that $\lambda_x$ is a monomorphism for every $x\in J$. Choose
a neighborhood $U_3$ of $(1,\mu)$ in $\Gl(V)\times {\mathcal M}$ such that if $(g',\mu')\in U_3$,
then $(g'\cdot\mu'-\mu)\in J$. Assume now that $(g,m)\in U_3$. Then we have $g\cdot m=\mu+a+b$
and hence $x=(a+b)\in J$. We therefore obtain by the above $\lambda_x(b)=0$. Since
$\lambda_x$ is a monomorphism, it follows that $b=0$. Thus $g\cdot m=\mu+a$ and
$\tau(g\cdot m)=\tau(\mu)$ as claimed.

\noindent{\bf (5)} Now we put together all ingredients in order to prove the theorem.
By step {\bf (3)}, $d_{(1,\mu)}\pi$ is surjective, and we may apply
Theorem \ref{inverse_function_theorem} again to obtain that for any neighborhood $U$ of
$(1,\mu)$ in $N'$, the image $\pi(U)$ is a neighborhood of $\pi(1,\mu)$.
$N'\subset\beta^{-1}(\tau(\mu))$ means that elements $(g,m)$ of $N'$ satisfy
$\pi_Z(\tau(g\cdot m))=\tau(\mu)$. By step {\bf (4)},
we can therefore suppose that all the elements $(g,m)$ of $U$ satisfy $\tau(g\cdot m)=\tau(\mu)$.
This is the claim of the stability theorem.
\end{prf}


\section{A rigid Leibniz algebra with non-trivial $HL^2$}  \label{section_example}


In this section, we finally obtain an analogue of Richardson's theorem \cite{Ric}
for Leibniz algebras.

Let ${\mathfrak g}$ be a semisimple Lie algebra and $M$ be an irreducible
(left) ${\mathfrak g}$-module (of dimension $\geq 2$). We put
$\widehat{\mathfrak g}=\mathfrak g\ltimes M$ semidirect product of $\mathfrak g$ and $M$
(with $[M,M]=0$).

\begin{theo}[Richardson's Theorem]  \label{Richardsons_theorem}
Let $\widehat{\mathfrak g}=\mathfrak g\ltimes M$ as above. Then $\widehat{\mathfrak g}$
is not rigid if and only if there exists a semisimple Lie algebra $\widehat{\mathfrak g}'$
which satisfies the following conditions:
\begin{enumerate}
\item[(a)] there exists a semisimple subalgebra ${\mathfrak g}'$ of
$\widehat{\mathfrak g}'$ which is isomorphic to ${\mathfrak g}$,
\item[(b)] if we identify ${\mathfrak g}'$ with ${\mathfrak g}$, then
$\widehat{\mathfrak g}'\,/\,{\mathfrak g}'$ is isomorphic to $M$ as a
${\mathfrak g}$-module.
\end{enumerate}
\end{theo}

\begin{prf} This is Theorem 2.1 in \cite{Ric}.\end{prf}

Richardson shows in \cite{Ric} that for ${\mathfrak g}={\mathfrak s}{\mathfrak l}_2(\C)$
and for $M=M_k$ the standard irreducible ${\mathfrak s}{\mathfrak l}_2(\C)$-module of
dimension $k+1$ and highest weight $k$,
the Lie algebra $\widehat{\mathfrak g}$ is not rigid if and only if $k=2,4,6,10$.
The semisimple Lie
algebras $\widehat{\mathfrak g}'$ are in this case the standard semisimple Lie
algebras of dimension $6,8,10$ and $14$. It turns out that $\widehat{\mathfrak g}$ has
necessarily rank 2 (see \cite{Ric}), and these are all rank 2 semi-simple Lie algebras
($A_1\times A_1$, $A_2$, $B_2$ and $G_2$).

We will extend Richardson's theorem to Leibniz algebras in the following sense.
First of all, we will restrict to simple Lie algebras ${\mathfrak g}$.
Let $I$ be another irreducible (right) ${\mathfrak g}$-module
(of dimension $\geq 2$).
We also set ${\mathfrak h}=I\dot{+}\widehat{\mathfrak g}$ the
hemisemidirect product with the $\mathfrak g$-module $I$ (in particular
$[{\mathfrak g},I]=0$,
$[I,{\mathfrak g}]=I$ and $[M,I]=0$). Observe that ${\mathfrak h}$ is a non-Lie
Leibniz algebra with ideal of squares $I$ and with quotient Lie algebra
$\widehat{\mathfrak g}$.

So, we have ${\mathfrak h}={\mathfrak g}\oplus M\oplus I$ as vector spaces.
In all the following, we fix the complex vector space ${\mathfrak g}\oplus M\oplus I$
and we will be
interested in the different Leibniz algebra structures on ${\mathfrak g}\oplus M \oplus I$.

\begin{theo} \label{rigidity_theorem}
Let ${\mathfrak h}=(V,\mu)$ be a finite-dimensional complex Leibniz algebra, and let
${\mathfrak g}$ be a subalgebra of ${\mathfrak h}$ such that $E^2({\mathfrak h},{\mathfrak g},
V)=0$.

Then the Leibniz algebra ${\mathfrak h}$ is not rigid if and only if there
exists a Leibniz algebra ${\mathfrak h}'$ which satisfies the following
conditions:
\begin{enumerate}
\item[a)] There exists a semisimple Lie subalgebra ${\mathfrak k}$ of
${\mathfrak h}'$ with a semisimple Lie subalgebra ${\mathfrak g}'\subset{\mathfrak k}$
which is isomorphic to ${\mathfrak g}$, i.e. ${\mathfrak g}\cong{\mathfrak g}'$;

\item[b)] if we identify the subalgebra ${\mathfrak g}'$ with ${\mathfrak g}$ by this
isomorphism, then ${\mathfrak k}/{\mathfrak g}$ is isomorphic to $M$ as a
${\mathfrak g}$-module and ${\mathfrak h}'/{\mathfrak k}$ is isomorphic to $I$ as a
Leibniz antisymmetric ${\mathfrak g}$-module.
\end{enumerate}
\end{theo}

\begin{prf} Let $V:={\mathfrak g}\oplus M\oplus I$. We will consider Leibniz algebras
${\mathfrak h}'$ on this fixed vector space $V$, i.e. ${\mathfrak h}=(V,\mu)$ and
${\mathfrak h}'=(V, \mu')$.

{\bf ``$\mathbf{\Leftarrow}$''}
Suppose that there exists a Leibniz algebra ${\mathfrak h}'$
which satisfies the conditions $a)$ and $b)$.

We may assume that $$\mu(x, x')=\mu'(x, x'), \ \mu(x,m)=\mu'(x,m), \
\mu(m, x)=\mu'(m, x), \ \mu(i, x)=\mu'(i, x) $$
for all $x,x'\in {\mathfrak g}$, all $m\in M$ and all $i\in I$.

Putting $g_t(x)=x, \ g_t(m)=tm, \ g_t(i)=ti$, we have that $g_t\in{\rm Gl}(V)$ for
all $t\not=0$ and
$$\lim_{t\rightarrow 0}g_t\cdot\mu'=\mu.$$
Therefore, $L$ is not rigid.

{\bf ``$\mathbf{\Rightarrow}$''}
Let $Leib$ be the set of all Leibniz algebras defined on the vector space $V$.
We are considering the Leibniz subalgebra ${\mathfrak s}\,:=\,{\mathfrak g}$
of ${\mathfrak h}$. By Proposition \ref{computations}, it satisfies
the cohomological condition in order to apply the stability theorem.
From Theorem \ref{stability_theorem}, we therefore get the existence of a
neighborhood $U$ of $\mu$
in $Leib$ such that if $\mu_1\in U$, the Leibniz algebra
${\mathfrak h}_1=(V,\mu_1)$ is isomorphic
to a Leibniz algebra ${\mathfrak h}'=(V,\mu')$ which satisfies the following
conditions:
\begin{enumerate}
\item[1)] $\mu(x, x')=\mu'(x, x'),$
\item[2)] $\mu(x, m)=\mu'(x,m), \ \mu(m, x)=\mu'(m, x),$
\item[3)] $\mu(i, x)=\mu'(i, x)=0, \ \mu(x,i)=\mu'(x,i),$
\end{enumerate}
for all $x,x'\in {\mathfrak g}$, all $m\in M$ and all $i\in I$.

As ${\mathfrak h}$ is supposed to be non rigid, we may assume that
${\mathfrak h}'$ is not
isomorphic to ${\mathfrak h}$. Since ${\mathfrak h}$ is a non-Lie Leibniz
algebra, the Leibniz algebra ${\mathfrak h}'$ is
also a non-Lie algebra.
Therefore, the ideal of squares of the algebra ${\mathfrak h}'$
is also non zero. We denote it by $I'$.

From conditions 1) and 2) we conclude that $I'\cap ({\mathfrak g}+M)=\{0\}$, and
thus $J:=I'\cap I\neq\{0\}.$
The condition 3) implies that $J$ is left module over ${\mathfrak g}$.
Since $I$ is an irreducible left ${\mathfrak g}$-module,
we have $J=I$ and thus $I'=I$ as a vector spaces.

We conclude that passing to the quotient with respect to $I=I'$ places us
in the same situation as in the proof of Richardson's theorem \cite{Ric}
p. 341, line 10: One may compare the two quotient Lie algebra structures
$L$ and $L'$ on the vector space ${\mathfrak g}\oplus M$, and their radicals.
The radical $R$ of $L'$ is (as in Richardson's proof) either zero or not,
which leads either to a semi-simple Lie algebra $L'$ satisfying the wanted
conditions a) and b) thanks to the conditions 1), 2) and 3) above,
or to a contradiction (to the fact that ${\mathfrak h}'$ is not
isomorphic to ${\mathfrak h}$). Therefore the
non-rigidity of ${\mathfrak h}$ implies the existence of a Leibniz algebra
${\mathfrak h}'$ satisfying conditions a) and b).    
\end{prf}

\begin{cor}
The Leibniz algebra ${\mathfrak h}:=I_l\dot{+}(M_k\rtimes{\mathfrak s}{\mathfrak l}_2(\C))$
for two standard irreducible left ${\mathfrak s}{\mathfrak l}_2(\C)$-modules $M_k$ and $I_l$
of highest weights $k=2n$ and $l$ respectively with odd integer $n>5$ and odd $l>2$ is rigid and satisfies
$HL^2({\mathfrak h},{\mathfrak h})\not= 0$.
\end{cor}

\begin{prf}
As discussed earlier, Richardson shows in \cite{Ric} that
for ${\mathfrak g}={\mathfrak s}{\mathfrak l}_2(\C)$
and for $M_k$ the standard irreducible ${\mathfrak s}{\mathfrak l}_2(\C)$-module of
dimension $k+1$ and highest weight $k$,
the Lie algebra $\widehat{\mathfrak g}=M_k\rtimes{\mathfrak g}$ is not rigid if and only if $k=2,4,6,10$.
He also shows that in case the half-highest-weight $\frac{k}{2}=:n$ is an odd integer $n>5$,
the Lie algebra cohomology of the Lie
algebra $\widehat{\mathfrak g}$
is not zero. As candidate for our Leibniz algebra ${\mathfrak h}$ satisfying the claim of the
corollary,
we take as before ${\mathfrak h}=I_l\dot{+}\widehat{\mathfrak g}$ for some irreducible
${\mathfrak s}{\mathfrak l}_2(\C)$-modules $I_l\not={\mathfrak s}{\mathfrak l}_2(\C)$
and $M_k$ such that the half-highest-weight $n:=\frac{k}{2}>5$ is an odd integer and $l>3$ is odd.
By Theorem \ref{rigidity_theorem},
${\mathfrak h}$ is then rigid.

On the other hand, by Proposition \ref{cohomology_double_semi_direct},
$H^2(\widehat{\mathfrak g},\widehat{\mathfrak g})$ injects into $HL^2({\mathfrak h},{\mathfrak h})$, thus we obtain that
this Leibniz cohomology space is not zero.
\end{prf}


\begin{appendix}
\section{Appendix: Lie-rigidity versus Leibniz-rigidity}

We record in this appendix further results on the question whether $H^2({\mathfrak g},{\mathfrak g})=0$ implies $HL^2({\mathfrak g},{\mathfrak g})=0$. The base field is fixed to be the field $\C$ of complex numbers. Recall that {\it rigid} stands for geometrically rigid, while we express the condition of being algebraically rigid by stating explicitly $H^2({\mathfrak g},{\mathfrak g})=0$. Algebraic rigidity implies geometric rigidity, but the converse is not true in general. Geometric rigidity means that any deformation must be isomorphic to the Lie algebra we started with.    

We have seen in Corollary \ref{proposition_bakhrom} that all (non nilpotent) solvable Lie algebras ${\mathfrak g}$ with $H^2({\mathfrak g},{\mathfrak g})=0$ and $\dim\,{\mathfrak q}>1$ satisfy
$HL^2({\mathfrak g},{\mathfrak g})=0$. In case the Lie algebra ${\mathfrak g}$ is only (Lie-)rigid, but does not necessarily satisfy $H^2({\mathfrak g},{\mathfrak g})=0$, we cannot conclude that
$HL^2({\mathfrak g},{\mathfrak g})=0$. But we will see in Theorems \ref{thm1} and \ref{thm2} below that for Lie algebras of a special form, one can still conclude in this situation that $Z({\mathfrak g})=0$.
In these special forms, we always suppose that the different pieces in the decomposition are non-zero. 

Recall that by Carles \cite{Car}, any rigid Lie algebra ${\mathfrak g}$ algebraic, i.e., it is isomorphic to the Lie algebra of an algebraic group. As the algebraicity implies the decomposability of the algebra \cite{Car}, it follows that for a rigid Lie algebra ${\mathfrak g}$, we have a decomposition ${\mathfrak g} = {\mathfrak s}\ltimes{\mathfrak r}$, where ${\mathfrak s}$ is a Levi subalgebra, ${\mathfrak r}={\mathfrak n}\rtimes {\mathfrak q}$ is the solvable radical of
${\mathfrak g}$, ${\mathfrak n}$ is the nilradical and ${\mathfrak q}$ is an exterior torus of derivations in the sense of Malcev; that is, ${\mathfrak q}$ is an abelian subalgebra of ${\mathfrak g}$ such that
${\rm ad}(x)$ is semisimple for all $x\in {\mathfrak q}$. Note that over the complex numbers, the semisimplicity of ${\rm ad}(x)$ means that there exists a basis of ${\mathfrak n}$ such that for any $x\in {\mathfrak q}$, the operator ${\rm ad}(x)_{|_{\mathfrak n}}$ has diagonal form and ${\rm ad}(x)_{|_{\mathfrak q}}=0$.

\begin{rem} \label{rem1} Since in the decomposition ${\mathfrak g} = {\mathfrak s}\ltimes{\mathfrak r}$ we have $[{\mathfrak s},{\mathfrak r}]\subseteq {\mathfrak r}$, we can view ${\mathfrak r}$ as an ${\mathfrak s}$-module. The fact that ${\mathfrak s}$ is semisimple implies that ${\mathfrak r}$ can be decomposed into a direct sum of irreducible submodules. Let ${\mathfrak r}_0$ be the sum of all $1$-dimensional submodules (so these are trivial submodules) and ${\mathfrak r}_{1}$ be the sum of the non-trivial irreducible submodules. Then
$$[{\mathfrak s},{\mathfrak r}_0]=0, \quad [{\mathfrak s},{\mathfrak r}_1]={\mathfrak r}_1, \quad [{\mathfrak r}_0,{\mathfrak r}_0]\subseteq {\mathfrak r}_0, \quad {\mathfrak r}_1 \subseteq {\mathfrak n}, \quad {\mathfrak q} \subseteq {\mathfrak r}_0.$$
\end{rem}

\begin{rem} Above, we have already considered the direct sum Lie algebra ${\mathfrak g}={\mathfrak s}{\mathfrak l}_2(\C)\oplus \C$. A short computation with the Hochschild-Serre spectral sequence shows that $H^2({\mathfrak g},{\mathfrak g})=0$, but (cf proof of Cor. 3 in \cite{FMM})
$\dim\,HL^2({\mathfrak g},{\mathfrak g})=1$.
Therefore, the Lie algebra ${\mathfrak g}$ is an example of a complex rigid Lie algebra of the form
${\mathfrak g}={\mathfrak s}\ltimes{\mathfrak n}$ where
$\dim\,({\mathfrak g}/[{\mathfrak g},{\mathfrak g}]) =1$ and where the center $Z({\mathfrak g})$ is
non-trivial. According to Theorem \ref{theorem_FMM}, the center is always non-trivial in case
$H^2({\mathfrak g},{\mathfrak g})$ and $HL^2({\mathfrak g},{\mathfrak g})$ differ.
\end{rem}

\begin{prop} \label{prop1} Let ${\mathfrak r}={\mathfrak r}_1\oplus \mathbb{C}^k$ be a split solvable Lie algebra. Then ${\mathfrak r}$ is not rigid.
\end{prop}

\begin{proof} Since ${\mathfrak r}$ is solvable, ${\mathfrak r}_1$ is also solvable and for this algebra we have ${\mathfrak r}_1={\mathfrak n}_1\oplus {\mathfrak q}_1$, where ${\mathfrak n}_1$ is the nilradical of ${\mathfrak r}_1$ and ${\mathfrak q}_1$ is a supplementary subspace. Let us fix elements $x\in {\mathfrak q}_1$ and $c=(c_1, c_2, \dots, c_k)\in \mathbb{C}^k$ and set $\varphi(x,c_i)=-\varphi(c_i,x)=c_i, 1\leq i \leq k$. Clearly, $\varphi \in Z^2({\mathfrak r},{\mathfrak r})$. Consider the infinitesimal deformation ${\mathfrak r}_t={\mathfrak r}+t\varphi$ of the algebra ${\mathfrak r}$ and let us show that $\varphi \not \in B^2(R,R)$.

Consider
$$(df)(x,c_i)=[f(x),c_i]+[x,f(c_i)]-f([x,c_i])=[x,f(c_i)], \ \ 1\leq i \leq k.$$

This equation implies that $\forall  f \in C^1({\mathfrak r},{\mathfrak r})$, we have $(df)(x,c_i)\neq c_i$, because
$[{\mathfrak r},{\mathfrak r}]\subset{\mathfrak n}_1$. Therefore, $\varphi \not \in B^2({\mathfrak r},{\mathfrak r})$ and ${\mathfrak r}_t={\mathfrak r}+t\varphi$ is a non-trivial deformation of the Lie algebra ${\mathfrak r}$. Therefore the Lie algebra ${\mathfrak r}$ is not rigid.
\end{proof}

\begin{cor} Let ${\mathfrak g}={\mathfrak s}\ltimes({\mathfrak n}\oplus {\mathfrak q})\oplus \mathbb{C}^k$ be a Lie algebra. Then ${\mathfrak g}$ is not rigid.
\end{cor}
\begin{proof} Due to Remark \ref{rem1}, we have that ${\mathfrak q}\cap[{\mathfrak g},{\mathfrak g}]=0$. So, taking ${\mathfrak g}_t={\mathfrak g}+t\varphi$ with $\varphi(x,c_i)=-\varphi(c_i,x)=c_i, 1\leq i \leq k$ similarly as in Proposition \ref{prop1}, we obtain a non-trivial deformation of ${\mathfrak g}$. Therefore, ${\mathfrak g}$ is not rigid.
\end{proof}

Let ${\mathfrak r}={\mathfrak n}\oplus{\mathfrak q}$ be a rigid solvable Lie algebra with basis
$\{x_1, x_2, \dots, x_k\}$ of ${\mathfrak q}$ and basis $N:=\{e_1, e_2, \dots, e_n\}$ of ${\mathfrak n}$.
Due to the arguments above, we can assume that its the table of multiplication has the following form:
$$\left\{\begin{array}{ll}
[e_i,e_j]=\sum\limits_{t=1}^{n}\gamma_{i,j}^te_{t},& 1\leq i,j\leq n,\\[1mm]
[e_i,x_j]=\alpha_{i,j}e_{i},& 1\leq i\leq n,\ \ 1\leq j\leq k, \ 1\leq k.\\[1mm]
\end{array}\right.$$
Note that for any $j\in \{1, \dots, k\}$, there exists $i$ such that $\alpha_{i,j}\neq 0$.

In a first part, we consider the center of rigid solvable Lie algebras of the form ${\mathfrak r}={\mathfrak n}\oplus{\mathfrak q}$. This study will cumulate in Theorem \ref{thm1} below.

\begin{prop}\label{prop2} Let ${\mathfrak r}={\mathfrak n}\oplus{\mathfrak q}$ be a rigid solvable Lie algebra such that $\dim\, {\mathfrak q}>1$. Then $Z({\mathfrak r})=\{0\}$.

\end{prop}
\begin{proof} First of all, we divide the basis $N$ of ${\mathfrak n}$ into two subsets. Namely,
$$N_1=\{ e_{i}\in N  \ | \ \exists j\in\{1, \dots, k\} \ \mbox{such  that} \ \alpha_{i,j}\neq0\},$$
$$N_2=\{ e_{i}\in N  \ | \ \alpha_{i,j}=0 \ \ \forall j\in\{1, \dots, k\}\}.$$

For convenience, we shall denote the bases $N_1$ and $N_2$ as follows:
$$N_1=\{e_1, e_2, \ldots, e_{n_1}\}, \quad N_2=\{f_1, f_2, \ldots, f_{n_2}\}.$$
Evidently, $n=n_1+n_2$.

Now let us suppose that $Z({\mathfrak r})\neq \{0\}$. Then there exists an element $0\neq c\in Z({\mathfrak r})$.
We put $$c=\sum\limits_{t=1}^{n_1}\beta_te_t+\sum\limits_{t=1}^{n_2}c_tf_t+\sum\limits_{t=1}^{k}b_tx_t.$$

Since the right bracketing  with $c$, i.e. $R_c$, satisfies $R_c\equiv 0$ and $x_1, \dots, x_k$ act diagonally on ${\mathfrak n}$, we derive
$b_t=0,\ \ 1\leq t\leq k$.

\textbf{Case 1.} Let $N_2=\emptyset$. Then $n=n_1$.

Consider
$$0=[c,x_i]=[\sum\limits_{t=1}^{n_1}\beta_te_t,x_i]=\sum\limits_{t=1}^{n}\beta_t\alpha_{t,i}e_t.$$
From this we deduce that $\beta_t=0$ for any $1\leq t\leq n_1$, as all $\alpha_{t,i}$ are non-zero.  Therefore, $c=0$. So, we have a contradiction with $Z({\mathfrak r})\not=\{0\}$.

\textbf{Case 2.} Let $N_2\neq\emptyset$. We set $\varphi(x_1,x_2)=c$, which is possible as by our hypothesis $\dim\,{\mathfrak q}>1$. It is easy to see that $\varphi$ is a $2$-cocycle.

Now we consider the associated infinitesimal deformation of the solvable Lie algebra ${\mathfrak r}_t={\mathfrak r}+t\varphi$. Let us check that this deformation ${\mathfrak r}_t$ is not equivalent to ${\mathfrak r}$. In fact, it is enough to check that $\varphi$ is not 2-coboundary of ${\mathfrak r}$.

Consider
$$(df)(x_1,x_2)=[f(x_1),x_2]+[x_1,f(x_2)]-f([x_1,x_2])=[f(x_1),x_2]+[x_1,f(x_2)].$$

Since ${\rm ad}(x_1)$ and ${\rm ad}(x_2)$ are diagonal and $c\in Z({\mathfrak r})$, we conclude that
$(df)(x_1,x_2)\neq c$ for any $f \in C^1({\mathfrak r},{\mathfrak r})$. This means, $\varphi \not \in B^2({\mathfrak r},{\mathfrak r}).$

Thus, the deformation ${\mathfrak r}_t={\mathfrak r}+\varphi$ is a non-trivial deformation of the rigid algebra ${\mathfrak r}$. This is a contradiction to the assumption that $Z({\mathfrak r})\not=\{0\}$.
\end{proof}

\begin{rem} Note that the rigidity of the solvable Lie algebra ${\mathfrak g}={\mathfrak n}\oplus {\mathfrak q}$ does not always imply that $\dim\,{\mathfrak q} > 1$. Indeed, there are examples of rigid solvable Lie algebras with $\dim \,{\mathfrak q}=1$ (see \cite{GAB}).
\end{rem}

Recall the commutation relations and notations which we introduced earlier:
$[e_i,x_j]=\alpha_{i,j}e_{i}$ for all $1\leq i\leq n$, all $1\leq j\leq k$, and all $1\leq k$, and
$\{e_1,\ldots,e_s\}$ are chosen Lie algebra generators of ${\mathfrak n}$.
In the following proposition, we suppress the index $j$ in $\alpha_{i,j}$ as ${\mathfrak q}$ is one dimensional and we have thus always $j=1$.

\begin{prop} \label{prop3} Let ${\mathfrak r}={\mathfrak n}\oplus{\mathfrak q}$ be a rigid solvable Lie algebra such that $\dim\,{\mathfrak q}=1$. Then $\alpha_i\neq 0$ for any $s+1\leq i \leq n.$
\end{prop}
\begin{proof} Let $<x>={\mathfrak q}$. Without loss of generality, one can assume that the basis elements $\{e_{1}, e_{2}, \dots, e_{n}\}$ of ${\mathfrak n}$ are homogeneous products of the generator basis elements, that is, for any $e_j\in [{\mathfrak n},{\mathfrak n}]$, we have
$$e_j =[\dots [e_{j_1}, e_{j_2}], \dots e_{j_{t_j}}]  \quad  s+1 \leq j \leq n, \quad 1 \leq j_i \leq s.$$

Then because of
$$[e_j,x]=(\alpha_{j_1} + \alpha_{j_2} + \dots + \alpha_{j_{t_j}})e_j,$$
we obtain
$$\alpha_j =\alpha_{j_1} + \alpha_{j_2} + \dots + \alpha_{j_{t_j}}  \quad s+1 \leq j \leq n, \quad 1 \leq j_i \leq s.$$

To the algebra ${\mathfrak r}$ as above and its regular element $x$, we associate a linear system $S(x)$ of equations of $n-1$ variables $z_1,\dots,z_{n-1}$ consisting of equalities $z_i + z_j = z_k$ if and only if the vector $[e_i,e_j]$ contains $e_k$ with a non-zero coefficient.
It is clear that $\{\alpha_1, \alpha_2, \dots, \alpha_n\}$ is a one of the solutions of this system.

If the system $S(x)$ has a unique fundamental solution, then we can obtain that  $z_p = k_p z_{i_0}$ with $k_p > 0.$ Therefore, in this case all elements of the solution are strictly positive. Since
 $\{\alpha_1, \alpha_2, \dots, \alpha_n\}$ is a solution, $\alpha_i\neq 0$ for any $s+1\leq i \leq n$.

So it remains to consider the case of a system of equations which has a space of fundamental solutions which is of dimension at least two.

Without loss of generality, we can assume $\alpha_1 \neq 0$ and let us assume that there exists some $p\geq s+1$ such that $\alpha_p=0.$ It follows then that $$\alpha_{p} =\alpha_{p_1} + \alpha_{p_2} + \dots + \alpha_{p_{t_p}}=0.$$

Let $\{z_{1}, z_{2}, \dots, z_{h}\}$ be another solution of $S(x)$, linearly independent from
$\{\alpha_1, \alpha_2, \dots, \alpha_n\}$. As $\alpha_1\neq 0$, we can assume $z_1=0$.

Consider the following solution $\{0, 1, 1, \dots, 1, z_{h+1},  \dots, z_{s}, z_{s+1}, z_{n}\},$ that is, here we get
  $z_{1}=0, z_{2}=1, \dots, z_{h}=1$.

Then we consider the following cochain:
$$\varphi(e_2, x) = e_2, \quad  \varphi(e_3, x) = e_3,\quad   \dots, \quad  \varphi(e_h, x) = e_h,$$
$$\varphi(e_{h+1}, x) = z_{h+1}e_{h+1}, \quad \varphi(e_{h+2}, x) = z_{h+2}e_{h+2}, \quad  \dots, \quad  \varphi(e_s, x) = z_s e_s,$$
$$\varphi(e_{s+1}, x) = z_{s+1}e_{s+1}, \quad  \varphi(e_{s+2}, x) = z_{s+2}e_{s+2}, \quad  \dots, \quad  \varphi(e_n, x) = z_n e_n.$$

It is easy to check that 
$\varphi \in Z^2({\mathfrak r},{\mathfrak r})$.

However, $\varphi \notin B^2({\mathfrak r},{\mathfrak r}).$ Indeed, if we had $\varphi \in B^2({\mathfrak r},{\mathfrak r}),$ then $\varphi = df$ for some $f \in C^1({\mathfrak r},{\mathfrak r})$.

Consider
$$e_2 = \varphi(e_2, x) = f([e_2, x]) - [f(e_2), x] - [e_2, f(x)]=$$
$$f(\alpha_2 e_2) - [c_2x+ \sum_{i=1}^{n}a_{2,i}e_i,x] - [e_2, d_0x + \sum_{i=1}^{n}d_ie_i]=$$
$$=\alpha_2(c_2x+ \sum_{i=1}^{n}a_{2,i}e_i) - \sum_{i=1}^{n}\alpha_ia_{2,i}e_i+d_0\alpha_2 e_2 + (*)= d_0\alpha_2 e_2 + (**).$$
Hence, we obtain $d_0\alpha_2=1$

On the other hand, we have
$$0 = \varphi(e_1, x) = f([e_1, x]) - [f(e_1), x] - [e_1, f(x)]=$$
$$f(\alpha_1 e_1) - [c_1x+ \sum_{i=1}^{n}a_{1,i}e_i,x] - [e_1, d_0x + \sum_{i=1}^{n}d_ie_i]=$$
$$=\alpha_1(c_1x+ \sum_{i=1}^{n}a_{1,i}e_i) - \sum_{i=1}^{n}\alpha_ia_{1,i}e_i+d_0\alpha_1 e_1 + (*)= d_0\alpha_1 e_1 + (**).$$

Therefore, $d_0\alpha_1=0.$ Since $\alpha_1\neq 0,$ we have $d_0=0.$ This is a contradiction to the assumption that $\varphi \in B^2({\mathfrak r},{\mathfrak r}).$
Finally, the deformation ${\mathfrak r}_t={\mathfrak r}+t\varphi$ is a non-trivial deformation of ${\mathfrak r}$, which contradicts the rigidity of ${\mathfrak r}$.
Therefore, we obtain that $\alpha_i\neq 0$ for any $1\leq i \leq n.$
\end{proof}

As a synthesis of Propositions \ref{prop2} and \ref{prop3}, we have the following theorem.

\begin{theo} \label{thm1} Let ${\mathfrak r}$ be a rigid solvable Lie algebra of the form ${\mathfrak r}={\mathfrak n}\oplus {\mathfrak q}$ with ${\mathfrak q}\not=0$. Then $Z({\mathfrak r})=0.$
\end{theo}
\begin{proof} The assertion of the theorem for the case when $\dim\,{\mathfrak q}>1$ follows from Proposition \ref{prop2}. Consider therefore the case $\dim\,{\mathfrak q}=1$. Let us assume that
$Z({\mathfrak r})\neq 0$. Then there exists $0\neq c\in Z({\mathfrak r})$. If $c\in {\mathfrak n} \setminus [{\mathfrak n},{\mathfrak n}]$, then we obtain that ${\mathfrak r}={\mathfrak r}_1\oplus \mathbb{C}$ and due to Proposition \ref{prop1}, we conclude that ${\mathfrak r}$ is not rigid.

Let now $c\in [{\mathfrak n},{\mathfrak n}]$. We set
$$c=\sum_{i=s+1}^{n}\beta_ie_i.$$
Taking into account Proposition \ref{prop3} and the following equality
$$0=[c,x]= \sum_{i=s+1}^{n}\beta_i[e_i,x]=\sum_{i=s+1}^{n}\beta_i\alpha_ie_i,$$
we conclude that $c=0$. The proof is complete.
\end{proof}

{\bf Conclusion 1.} Let ${\mathfrak r}$ be a solvable Lie algebra of the form ${\mathfrak r}={\mathfrak n}\oplus{\mathfrak q}$ such that $H^2({\mathfrak r},{\mathfrak r})=0$, then ${\mathfrak r}$ is rigid. By Theorem \ref{thm1}, we have $Z({\mathfrak r})=\{0\}$. Now applying Theorem \ref{theorem_FMM},
we conclude that $HL^2({\mathfrak r},{\mathfrak r})=0$, which implies by  \cite{Bal} that ${\mathfrak r}$ is rigid as a Leibniz algebra.\\

In a second part, we will now study the center of general rigid Lie algebras of the form
${\mathfrak g}={\mathfrak s}\ltimes({\mathfrak n}\oplus{\mathfrak q})$ using the same methods as before.

\begin{prop}\label{prop4} Let ${\mathfrak g}={\mathfrak s}\ltimes({\mathfrak n}\oplus {\mathfrak q})$ be a rigid Lie algebra such that $\dim\,{\mathfrak q}>1$. Then $Z({\mathfrak g})=\{0\}$.
\end{prop}
\begin{proof} Let us suppose that $Z({\mathfrak g})\neq 0$ and $0\neq c\in Z({\mathfrak g})$. Clearly, $c\in {\mathfrak n}.$  From Remark \ref{rem1} we conclude that $[{\mathfrak g},{\mathfrak g}]\cap {\mathfrak q}=0.$ Then similarly as in the Proposition \ref{prop2} we conclude that
${\mathfrak g}_t={\mathfrak g}+t\varphi$ with $\varphi(x_1,x_2)=c$ is a non-trivial deformation of
${\mathfrak g}$. This implies that $Z({\mathfrak g})=\{0\}.$
\end{proof}

As recalled in the beginning, if ${\mathfrak g}={\mathfrak s}\ltimes({\mathfrak n}\oplus {\mathfrak q})$ is rigid, then ${\mathfrak s}$ is a Levi subalgebra of ${\mathfrak g}$, ${\mathfrak n}$ is the nilradical and ${\mathfrak q}$ is an abelian subalgebra with diagonal operators ${\rm ad}(x)_{|_{\mathfrak n}}$ for any $x\in{\mathfrak q}$. Moreover, $[\mathfrak{h},\mathfrak{h}]=[{\mathfrak q},{\mathfrak q}]=0$ for a Cartan subalgebra ${\mathfrak h}$ of ${\mathfrak s}$. Now due to Remark \ref{rem1} (indeed, since ${\mathfrak q}\subseteq {\mathfrak r}_0$, we obtain $[{\mathfrak s},{\mathfrak q}]\subseteq [{\mathfrak s},{\mathfrak r}_0]=0$), we conclude that ${\mathfrak h}\oplus {\mathfrak q}$ is toroidal subalgebra in ${\mathfrak g}$ which acts on
$${\mathfrak n}\oplus (\oplus_{i} {\mathfrak l}_{\beta_i})\oplus (\oplus_{i} {\mathfrak l}_{-\beta_i}),$$
where
$${\mathfrak s}=\mathfrak{h}\oplus (\oplus_{i} {\mathfrak l}_{\beta_i})\oplus (\oplus_{i} {\mathfrak l}_{-\beta_i}).$$

Therefore, we have the following toroidal decomposition:
$${\mathfrak g}=(\mathfrak{h}\oplus {\mathfrak q})\oplus({\mathfrak g}_{\gamma_1}\oplus {\mathfrak g}_{\gamma_2} \oplus \dots \oplus {\mathfrak g}_{\gamma_t})$$
such that
$[g,h+x]=\gamma(h+x)g$ for any $h\in \mathfrak{h}, x\in {\mathfrak q}$ and $g\in {\mathfrak g}_\gamma$.

Let $\dim\,{\mathfrak q}=1$ with ${\mathfrak q}=<x>$, $\dim\,\mathfrak{h}=q$ and ${\mathfrak l}_{\beta_i}=<s_{\beta_i}>$ and ${\mathfrak l}_{-\beta_i}=<s_{-\beta_i}>$.
Then we have the following products in ${\mathfrak g}$:
$$\left\{\begin{array}{ll}
[e_i,e_j]=\sum\limits_{t=1}^{n}\gamma_{i,j}^te_{t},& 1\leq i,j\leq n,\\[1mm]
[e_i,x]=\alpha_{i}e_{i},& 1\leq i\leq n,\\[1mm]
[e_i,h_p]=\theta_{i,p}e_{i},& 1\leq i\leq n,\ \ 1\leq p\leq q,\\[1mm]
[x,{\mathfrak s}]=0.& \\[1mm]
\end{array}\right.
$$

\begin{prop}\label{prop5} Let ${\mathfrak g}={\mathfrak s}\ltimes({\mathfrak n}\oplus {\mathfrak q})$ be a rigid Lie algebra such that $\dim\,{\mathfrak q}=1$. Then $Z({\mathfrak g})=\{0\}$.
\end{prop}
\begin{proof}
Similarly as in the proof of Proposition \ref{prop3} we consider the system of linear equations with respect to $\alpha_{1}, \alpha_{2}, \dots, \alpha_{s}:$
$$\alpha_j =\alpha_{j_1} + \alpha_{j_2} + \dots + \alpha_{j_{t_j}}  \quad s+1 \leq j \leq n, \quad 1 \leq j_i \leq s.$$

Let us assume that there exists some $p\geq s+1$ such that $\alpha_p=0.$ It follows that $$\alpha_{p} =\alpha_{p_1} + \alpha_{p_2} + \dots + \alpha_{p_{t_p}}=0.$$

To the above relations, we associate the system $S(z)$ of linear equations with respect to $z_{1}, z_{2}, \dots, z_{s}$. That is, the system contains the relations
$$z_j =z_{j_1} + z_{j_2} + \dots + z_{j_{t_j}}  \quad s+1 \leq j \leq n, \quad 1 \leq j_i \leq s.$$

Note that the fundamental solution of the system $S(z)$ does not consist of a unique solution. Indeed, if it had a
unique fundamental solution, then we would obtain that $0=z_p = k_p z_{1}$ with $k_p \neq 0$, which implies $z_1=0$
and hence, $\alpha_{1} = \alpha_{2} = \dots = \alpha_{n}=0$ which is impossible because in this case ${\rm ad}(x)$ acts trivially to ${\mathfrak n}$.

So, we conclude that the system $S(z)$ of equations has a fundamental solution which is a vector space of dimension at least two.
Let us suppose that $\{z_{1}, z_{2}, \dots, z_{h}\}$ is another solution of $S(z)$, linearly independent of $\{\alpha_1,\ldots,\alpha_h\}$. Moreover, we may assume that $\alpha_{1}$ appears in the set $\{\alpha_{p_1}, \alpha_{p_2}, \dots,\alpha_{p_{t_p}}\},$ i.e.,
$e_{1}$ appears in the following product
$$e_p = [\dots [e_{p_1}, e_{p_2}], \dots e_{p_{t_p}}]$$

Consider following solution $\{0, 1, 1, \dots, 1, z_{h+1},  \dots, z_{s}, z_{s+1},\dots, z_{n}\},$ that is, here we get
$z_{1}=0, z_{2}=1, \dots, z_{h}=1$.

We have
$$[[e_{i}, s_\beta],x]=[[e_{i}, x],s_\beta]+[e_{i}, [s_\beta,x]]=[[e_{i}, x],s_\beta]=\alpha_i[e_{i}, s_\beta].$$

Then we consider the following cochain:
$$\varphi(e_2, x) = e_2, \quad  \varphi(e_3, x) = e_3,\quad   \dots, \quad  \varphi(e_h, x) = e_h,$$
$$\varphi(e_{h+1}, x) = z_{h+1}e_{h+1}, \quad \varphi(e_{h+2}, x) = z_{h+2}e_{h+2}, \quad  \dots, \quad  \varphi(e_s, x) = z_s e_s,$$
$$\varphi(e_{s+1}, x) = z_{s+1}e_{s+1}, \quad  \varphi(e_{s+2}, x) = z_{s+2}e_{s+2}, \quad  \dots, \quad  \varphi(e_n, x) = z_n e_n.$$
$$\varphi(e_{i}, h_p) = \theta_{i,p}e_{i}, \quad \varphi(h_p,x)=0, \quad  1\leq i\leq n,\ 1\leq p\leq q,$$
$$\varphi(e_{i}, s_\beta) =\varphi(x,s_\beta)=\varphi(e_{i}, e_{j})=0, \quad  \varphi([e_{i}, s_\beta],x)=z_i[e_{i}, s_\beta], \quad 1\leq i,j\leq n,\ 1\leq p\leq q.$$


By straightforward computations we check that $\varphi$ is a $2$-cocycle of ${\mathfrak g}$.
%
%
%
%
%
%


However, $\varphi \notin B^2({\mathfrak g},{\mathfrak g}).$ Indeed, if $\varphi \in B^2({\mathfrak g},{\mathfrak g}),$ then $\varphi = df$ for some $f=f_{\mathfrak r}+f_{\mathfrak s}\in C^1({\mathfrak g},{\mathfrak g})$ with $f_{\mathfrak r}: {\mathfrak g}\rightarrow {\mathfrak r}$ and $f_{\mathfrak s}: {\mathfrak g}\rightarrow {\mathfrak s}.$

Consider
$$0=\varphi({\mathfrak s}, x) = f([{\mathfrak s}, x]) - [f({\mathfrak s}), x] - [{\mathfrak s},f(x)]=-[f({\mathfrak s}),x]-[{\mathfrak s},f_{\mathfrak r}(x)]-[{\mathfrak s},f_{\mathfrak s}(x)].$$
This implies $[f({\mathfrak s}),x]+[{\mathfrak s},f_{\mathfrak r}(x)]=0$ and $[{\mathfrak s},f_{\mathfrak s}(x)]=0$.

Taking into account that for any $s\in {\mathfrak s}$, there exists $s'\in {\mathfrak s}$ such that
$[s,s']\neq 0$, we conclude that $f_{{\mathfrak s}}(x)=0.$

Consider
$$e_2 = \varphi(e_2, x) = f([e_2, x]) - [f(e_2), x] - [e_2, f(x)]=$$
$$f(\alpha_2 e_2) - [c_2x+ \sum_{i=1}^{n}a_{2,i}e_i+f_{S}(e_2),x] - [e_2, d_0x + \sum_{i=1}^{n}d_ie_i]=$$
$$=\alpha_2(c_2x+ \sum_{i=1}^{n}a_{2,i}e_i+f_{S}(e_2)) - \sum_{i=1}^{n}\alpha_ia_{2,i}e_i+d_0\alpha_2 e_2 + (*)= d_0\alpha_2 e_2 + (**).$$
Hence, we obtain $d_0\alpha_2=1$

On the other hand, we have
$$0 = \varphi(e_1, x) = f([e_1, x]) - [f(e_1), x] - [e_1, f(x)]=$$
$$f(\alpha_1 e_1) - [c_1x+ \sum_{i=1}^{n}a_{1,i}e_i+f_{S}(e_1),x] - [e_1, d_0x + \sum_{i=1}^{n}d_ie_i]=$$
$$=\alpha_1(c_1x+ \sum_{i=1}^{n}a_{1,i}e_i+f_{S}(e_1)) - \sum_{i=1}^{n}\alpha_ia_{1,i}e_i+d_0\alpha_1 e_1 + (*)= d_0\alpha_1 e_1 + (**).$$

Therefore, $d_0\alpha_1=0.$ Since $\alpha_1\neq 0,$ we have $d_0=0.$ This is a contradiction to assumption that $\varphi \in B^2({\mathfrak g},{\mathfrak g}).$

Finally, the deformation ${\mathfrak g}_t={\mathfrak g}+t\varphi$ is a non-trivial deformation of
${\mathfrak g}$, which contradicts the rigidity of ${\mathfrak g}$.

Thus, we have proved that if $[e_p,x]=0$ for some $e_p\in {\mathfrak n}$, then ${\mathfrak g}$ is non-rigid. Similar as in the proof of Theorem \ref{thm1}, we obtain that $Z({\mathfrak g})=\{0\}$.
\end{proof}

Again, we perform a synthesis of the preceding two propositions in the following theorem.

\begin{theo} \label{thm2} Let ${\mathfrak g}={\mathfrak s}\ltimes({\mathfrak n}\oplus{\mathfrak q})$ be a rigid Lie algebra such that ${\mathfrak q}\neq 0$. Then $Z({\mathfrak g})=0.$
\end{theo}
\begin{proof} The assertion of the theorem for the case when $\dim\,{\mathfrak q}>1$ follows from Proposition \ref{prop4}. In case $\dim\,{\mathfrak q}=1$, we have that $Z({\mathfrak g})=0$
due to Proposition \ref{prop5}.
\end{proof}

{\bf Conclusion 2.} Let ${\mathfrak g}={\mathfrak s}\ltimes({\mathfrak n}\oplus{\mathfrak q})$ be a Lie algebra such that $H^2({\mathfrak g},{\mathfrak g})=0$, then ${\mathfrak g}$ is rigid. From Theorem \ref{thm2} we obtain that $Z({\mathfrak g})=\{0\}$. Now applying Theorem \ref{theorem_FMM}, we conclude that $HL^2({\mathfrak g},{\mathfrak g})=0$, which implies by \cite{Bal} that ${\mathfrak g}$ is rigid as a Leibniz algebra.

\end{appendix}

\end{document}